# A Survey on Convex Optimization for Guidance and Control of Vehicular Systems

Zhenbo Wang[a,*]

[a]*Department of Mechanical, Aerospace, and Biomedical Engineering, The University of Tennessee, Knoxville, Tennessee, 37996, USA*

**Abstract**

Guidance and control (G&C) technologies play a central role in the development and operation of vehicular systems. The emergence of computational guidance and control (CG&C) and highly efficient numerical algorithms has opened up the great potential for solving complex constrained G&C problems onboard, enabling higher level of autonomous vehicle operations. In particular, convex-optimization-based G&C has matured significantly over the years and many advances continue to be made, allowing the generation of optimal G&C solutions in real-time for many vehicular systems in aerospace, automotive, and other domains. In this paper, we review recent major advances in convex optimization and convexification techniques for G&C of vehicular systems, focusing primarily on three important application fields: 1) Space vehicles for powered descent guidance, small body landing, rendezvous and proximity operations, orbital transfer, spacecraft reorientation, space robotics and manipulation, spacecraft formation flying, and station keeping; 2) Air vehicles including hypersonic/entry vehicles, missiles and projectiles, launch/ascent vehicles, and low-speed air vehicles; and 3) Motion control and powertrain control of ground vehicles. Throughout the paper, we draw figures that illustrate the basic mission concepts and scenarios and present tables that summarize similarities and distinctions among the key problems, ideas, and approaches. Where available, we provide comparative analyses and reveal correlations between different applications.

*Corresponding author

*Email address:* zwang124@utk.edu (Zhenbo Wang)

Finally, we identify open challenges and issues, discuss potential opportunities, and make suggestions for future research directions.

*Keywords:* Guidance and Control, Optimal Control, Trajectory Optimization, Convex Optimization, Space Vehicles, Air Vehicles, Ground Vehicles

## 1. Introduction

Guidance and control (G&C) technologies are vital for vehicular systems. In recent years, the generation of G&C commands relies much more extensively on onboard computation, accelerated by a critical need of highly efficient G&C systems for intelligent, autonomous vehicle operations. For example, NASA has been driving and supporting innovation in autonomous systems and key technologies such as flight computing and G&C for the development of future air transportation and new space exploration missions to the Moon, Mars, and other planetary bodies, as highlighted in the 2015 NASA Technology Roadmaps [1] and the more recent 2020 NASA Technology Taxonomy [2]. In January 2023, the United States Department of Defense (USDOD) updated its Directive 3000.09, Autonomy in Weapon Systems, which shows USDOD's commitment to developing, testing, fielding, and employing autonomous and semi-autonomous weapon systems [3]. Furthermore, the United States Department of Transportation (USDOT) has been devoted to enabling safe, efficient, and sustainable transportation systems by promoting new forms of mobility and advancing transportation technologies from electric, automated, and connected vehicles to advanced air mobility and commercial space travel [4, 5]. The overall goal of these efforts is to facilitate paths to enabling more efficient and capable space, air, and ground vehicle operations in various mission scenarios. The objectives are to transform next-generation vehicular systems from manually controlled systems to ones that respond to the dynamic mission requirements in real-time and operate autonomously in highly uncertain environments. Towards achieving this, an accurate and robust solution process needs to be created for the vehicle to perform complex missions incorporating highly nonlinear vehicle dynamic systems and stringent constraints with a high degree of reliability. This process will involve automating and optimizing some of the system functions such as G&C to enable fast yet accurate decision-making systems.

The field of G&C has recently been evolving from focusing on traditional laws and controllers to numerical algorithms with the aim of achieving onboard applications for autonomous vehicle systems [6]. Specifically, an emerging and accelerating trend has occurred in G&C, where the traditional algebraic G&C laws are replaced by numerical and computational algorithms. In contrast to traditional G&C, computational guidance and control (CG&C) allows complex G&C missions involving highly nonlinear dynamic systems and many state and control constraints to be performed. It is worth noting that CG&C is not simply solving G&C problems numerically onboard. Reliability, accuracy, computational efficiency, and robustness of the solution process are all primary challenges facing the G&C community in the development of numerical G&C algorithms. Another observation is that the vehicle G&C system has been generally decoupled, in which the guidance subsystem determines the desired trajectory as well as the associated changes in position, velocity, rotation, and acceleration for the vehicle to move from its current location to a designated target, while the control subsystem manipulates the forces acting on the vehicle via steering control defectors such as aerodynamic surfaces and thrusters to execute the guidance commands while maintaining vehicle stability [7]. Such an approach breaks the entire G&C problem down into a series of subproblems including reference trajectory generation, tracking guidance, attitude control, and iterative online implementation. The large majority of these problems have been characterized by and intrinsically tied with optimal control problems (OCPs). In many cases, unfortunately, analytical solutions are usually impossible to find, and numerical techniques must be employed to determine feasible reference trajectories and closed-loop control policies.

So far, the OCP-based G&C problems have been substantially addressed by either the indirect method or the direct method. The pros and cons of these methods have been extensively discussed in the literature, and details can be found in many publications such as [8] and [9]. Briefly speaking, the indirect method builds on the calculus of variations and Pontryagin's Minimum Principle, derives the neces-sary conditions such as adjoint equations and transversality conditions, determines the optimal solutions by minimizing the Hamiltonian with respect to the control, and reduces the initial OCP to a multi-point boundary value problem. The optimality of the indirect method can be guaranteed; however, complicated and lengthy mathematical derivations are needed, and high-quality initial guesses of the adjoints are always required. In contrast, the direct method does not require explicit derivation of optimality conditions; instead, it discretizes the continuous trajectory into multiple segments and converts the original continuous-time OCP into a finite-dimensional parameter optimization problem, which is then solved using numerical optimization methods such as nonlinear programming (NLP). The direct method is easy to implement; for complicated, highly nonconvex problems, however, the solution process is usually time-consuming, and the convergence of NLP algorithms is hard to be guaranteed. Therefore, despite decades of advancement, there is still a lack of highly efficient algorithms that are capable of handling highly nonlinear system dynamics and a variety of mission constraints with stable convergence and real-time performance without compromising solution accuracy and optimality for G&C across all mission phases in both single-vehicle and multi-vehicle settings. Convex-optimization-based G&C has emerged and advanced in the past two decades, providing great potential to address these issues and challenges and achieve the collective goals of autonomous vehicle systems. With the significant increase in computational efficiency, convex-optimization-based G&C is expected to become a fundamental technology for autonomous system operations.

As a subfield of mathematical optimization that addresses the problem of minimizing convex functions over convex sets, the study of convex optimization dates back to more than a century ago; however, the power of convex optimization for practical applications did not come to light until the 1990s, when it was discovered that many engineering problems are actually convex or can be approximated as convex optimization problems [10]. Recent development of highly efficient convex optimization algorithms together with the advances of computing power have provided the basis for substantial increase in the performance of generating optimal vehicle trajectories and producing closed-loop tracking control command due to its advantages including low complexity, polynomial-time computation, global optimality, and deterministic convergence [10, 11]. In addition to allowing more rapid and stable system operations, these algorithmic advances offer the potential for automating several of the G&C tasks, elevating the human pilot or driver to the role of operation managers, where the human is expected to intervene only if the automated system is unable to deal with the situation at hand. For these reasons, NASA cited convex optimization as a computationally-efficient method for solving large divert guidance problems in real-time for potential future entry, descent, and landing (EDL) applications [1]. Nevertheless, most real-world problems are nonconvex and difficult to solve in both theory and practice. While a nonconvex problem can be potentially addressed using techniques commonly known as convexification and relaxed into a hierarchy of convex subproblems that can be reliably solved using efficient interior-point methods (IPMs), the obtained candidate solutions (if converged) can only be suspected of being locally optimal solutions, and the optimality of the solutions is generally difficult to validate [12].

It is generally acknowledged that the initial impetus for the development of convex optimization algorithms for vehicle G&C applications began in the U.S. with a series of publications on powered descent guidance for Mars pinpoint landing [13, 14]. In the aerospace domain, convex optimization was initially used to solve the optimal powered descent guidance problem, where a propellant-optimal trajectory optimization problem was formulated as an OCP subject to state and control constraints and relaxed into a convex optimization problem through the lossless convexification technique [15, 16, 17]. It was then followed by the development of basic convex optimization and sequential convex programming (SCP) algorithms and their more advanced variants with enhanced techniques such as virtual control and pseudospectral discretization for more applications including rendezvous and proximity operations [18, 19], low-thrust orbital transfers [20, 21], and space robotics [22]. Convex-optimization-based G&C algorithms for high-speed atmospheric flight vehicles appeared around 2015 [23, 24], and the past five years observed rapid advances in algorithm development with a primary growth of SCP-type algorithms and their applications to hypersonic/entry vehicles [25, 26, 27, 28] and launch/ascent ve-hicles [29, 30, 31, 32]. Convex optimization for low-speed air vehicles appears to be following a similar developmental pathway. Research in earlier years between 2000 and 2010 was mainly focused on solving single convex optimization problems with simplified vehicle models to obtain the solution [33, 34], while SCP is at the heart of the development of more advanced algorithms for wide air vehicle applications in recent years [35, 36, 37, 38]. More recently, convex optimization has gained significant interest in the automotive domain to improve the efficiency and performance of ground vehicles. The major focus is on the generation of approximate optimal solutions to motion/speed control [39, 40, 41] and powertrain control [42, 43] problems in the context of smart mobility and intelligent transportation systems.

The primary goal of this paper is to provide a holistic survey on the development of convex optimization algorithms for G&C of vehicular systems. In particular, this paper presents fundamental results and latest advances in convexification and SCP techniques for space, air, and ground vehicles, and highlights some limitations of the existing solutions in each of these areas from both theoretical and technological perspectives. Finally, this survey paper presents some relevant challenges and issues that hinder the implementation of convex-optimization-based methods for real-world G&C missions, and discuss potential research efforts to reduce these challenges and issues in the future. It is worth mentioning that this paper not only surveys the traditional G&C areas but also reviews the applications of convex optimization techniques for G&C in newly emerged fields such as reusable rocket landing, small body exploration, electric vertical take-off and landing (eVTOL) vehicles, advanced air mobility (AAM), connected and automated vehicles (CAVs), and collaborative space-air-ground missions. The intended readers for this paper are researchers, students, and professionals who are interested in the design of G&C systems for space, air, or ground vehicles and those who would like to explore the applications of convex optimization for other areas. Towards achieving these goals, this paper surveys the most representative works published in the past 20 years with preference to journal and international conference papers. The fundamental ideas underlying the most popular results are discussed in great detail, whereas the follow-on, incremental results are merely cited. Higher priority is given to the breadth rather than the depth of results presented. Also, skipping the discussions of the classical optimal control theory as well as the indirect and direct methods that have been extensively analyzed in numerous survey papers such as [8] and [9], this paper focuses on convex optimization and SCP methods, introduces the algorithms at a high level, and chooses to cover wider applications ranging from space to air and to ground vehicles. The motivation for covering a larger number of applications comes from the fact that methods developed to solve problems from different domains are closely interconnected and mutually inspired from both theoretical development and practical application perspectives.

This survey paper complements and extends the existing survey and review papers on guidance, control, and trajectory optimization. For example, [13] provided an overview of some common convexification techniques and their applications to aerospace G&C problems solved by early 2017. Over the past years, there was an exponential growth of publications on convex optimization for wider G&C application domains. In fact, over 300 publications have been added to the convex-optimization-based G&C literature in the past seven years and have significantly advanced this area such that new theoretical results need to be synthesized and new applications need to be assessed. Also, [14] surveyed the general optimization-based methods for space vehicle control with a focus on the last ten years of advances in convex optimization techniques for G&C of space vehicles including launchers, planetary landers, satellites, and spacecraft. However, this paper did not cover some related topics such as G&C of purely atmospheric vehicles (e.g., missiles and hypersonic aircraft), satellite swarms, low-speed air vehicles, novel mobility concepts, the emerging CAVs, and so on. The interconnections among these areas need to be explored and further challenges and issues need to be addressed. Moreover, [44] provided a comprehensive tutorial of lossless convexification, successive convexification, and guaranteed sequential trajectory optimization methods for reliable and efficient trajectory generation, accompanied by an open-source SCP toolbox as well as a number of case studies and numerical examples. In addition, interested readers are referred to [8] and [9] for surveys on numerical techniques for solving trajectory optimization problems and general OCPs via indirect and direct meth-ods; [45] for an introductory tutorial that covers the basics for numerically solving trajectory optimization problems with a focus on direct collocation methods; [46] for a survey on common transcription methods that convert the continuous OCP into a parameter optimization problem as well as the evolutionary algorithms or metaheuristics that can solve the resulting parameter optimization problem; [47] for a review of theoretical foundations of the pseudospectral optimal control theory for practical implementation in aerospace and autonomous systems; [48] for a survey on mathematical techniques, including geometric optimal control, continuation/homotopy method, and dynamical system theory, to improve the performance of optimal control tools such as the Pontryagin's Minimum Principle in solving OCPs in aerospace; [49] for a survey on the guidance methods, including analytical guidance methods, numerical optimization algorithms, convexification strategies, and learning-based methods, for pinpoint soft-landing on the Moon, Mars, and Earth; [50] and [51] for surveys on the state-of-the-art in G&C techniques for interception, descent, and landing on small celestial bodies such as asteroids and comets; [52] for a survey on optimization approaches to civil applications of unmanned aerial vehicles (UAVs); [53] for a survey on planning and control of CAVs with a particular focus on approaches to improving energy efficiency; [54] and [55] for surveys on model predictive control (MPC) for aerospace systems and general dynamical systems, respectively; and [12] for a survey on major advances in conic optimization and its applications in machine learning, power systems, state estimation, and the abstract problems of rank minimization and quadratic optimization. Finally, readers interested in detailed results are referred to the literature reviews presented in each of the papers cited herein.

This survey paper is organized as follows. Section 2 presents a brief overview of the general OCP formulation for G&C applications and the basic convex optimization and SCP algorithms as well as some enhanced techniques for solving these problems. Section 3 surveys the applications of convex optimization for G&C of space vehicles. This section is structured by considering powered descent guidance first and then rendezvous and proximity operations, orbital transfers, and so on. Section 4 surveys the applications of convex optimization for G&C of air vehicles considering both high-speed and low-speed vehicles under multiple atmospheric flight missions and scenarios. Section 5 gives special emphasis to the application of convex optimization for G&C of ground vehicles and highlights some convex-optimization-based techniques used to improve ground mobility efficiency. Section 6 discusses some open research challenges and issues and recommends some future research directions. Finally, Section 7 concludes this survey paper.

## 2. Convex Optimization and Sequential Convex Programming

Many vehicular G&C problems are formulated as optimal control problems (OCPs). This paper focuses on the survey of OCP-based G&C problems that can be solved by convex optimization algorithms. This section gives a brief introduction on OCP as well as the basic convex optimization and sequential convex programming (SCP) techniques.

### *2.1. Optimal Control Problem*

A continuous-time OCP can be generally posed as [56]:

*Problem 1*:

$$\text{minimize: } J = \Phi[\mathbf{x}(t_0), t_0, \mathbf{x}(t_f), t_f] + \int_{t_0}^{t_f} L[\mathbf{x}(t), \mathbf{u}(t), t]dt \tag{1}$$

$$\text{subject to: } \dot{\mathbf{x}}(t) = \mathbf{f}[\mathbf{x}(t), \mathbf{u}(t), t] \tag{2}$$

$$\boldsymbol{\varphi}_{\min} \le \boldsymbol{\varphi}[\mathbf{x}(t_0), t_0, \mathbf{x}(t_f), t_f] \le \boldsymbol{\varphi}_{\max} \tag{3}$$

$$\mathbf{C}_{\min} \le \mathbf{C}[\mathbf{x}(t), \mathbf{u}(t), t] \le \mathbf{C}_{\max} \tag{4}$$

where $t \in [t_0, t_f]$ is the independent variable that usually denotes time, $\mathbf{x}(t) \in \mathbf{R}^{n_x}$ is the state trajectory, and $\mathbf{u}(t) \in \mathbf{R}^{n_u}$ is the control history. Solving an OCP aims to determine the optimal control history $\mathbf{u}^*(t)$ that drives the system from an initial state $\mathbf{x}(t_0)$ at an initial time $t_0$ to a target state $\mathbf{x}$ at a terminal time $t_f$ while minimizing a performance measure Eq. (1) and satisfying the dynamics Eq. (2), boundary conditions Eq. (3), and path constraints Eq. (4). The initial state and initial time are usually specified, while the terminal state and terminal time can be free. The system dynamics in Eq. (2) can be described as a linear time-invariant, linear time-varying, nonlinear time-invariant, or nonlinear time-varying system, according to the specific form of state equations used. Also, simple bounds on the state and control variables are special cases of Eq. (4). Problem 1 can be cast as a minimum-time, terminal control, minimum-control-effort (e.g., minimum-fuel or minimum-energy), or tracking problem, depending on the form of the objective functional defined in Eq. (1).

The particular focus of this paper is convex-optimization-based G&C techniques within the scope of the direct optimal control method. There are many discretization methods that can transform the infinite-dimensional OCP in Problem 1 into a finite-dimensional numerical optimization problem [8, 9]. When a discretization method is employed, the continuous interval of the independent variable (usually time) is discretized, and the state and control histories are represented by sequences of discrete nodes. The system dynamics can be satisfied via explicit or implicit numerical integration and become equality constraints. All other constraints are also enforced at the discrete nodes. Eventually, the original OCP problem is transcribed into a parameter optimization problem that can be solved by an NLP solver. The readers are referred to [9], [46], [14], and many other publications for the detailed discretization process and methods. It is notable that the discretized problem may take different forms and fall into one of the sub-classes of convex optimization problems that can be solved very efficiently.

### *2.2. Convex Optimization*

If the discretized problem can be formulated as or relaxed into a convex optimization problem (called one-shot convexification as in Figure 1), such as linear programming (LP), convex quadratic programming (QP), convex quadratically constrained quadratic programming (QCQP), or second-order cone programming (SOCP), the problem can be solved in polynomial time because of its low complexity [10]. Then, state-of-the-art interior-point methods (IPMs) can be used to compute a globally optimal solution with deterministic stopping criteria and a prescribed level of accuracy [11]. Additionally, when solving a convex optimization problem, no initial guesses need to be supplied by users, because self-dual embedding techniques allow IPMs to start from a self-generated feasible point [57]. All of these characteristics offer the advantages not observed in the traditional direct or indirect method and provide great opportunities for onboard applications.

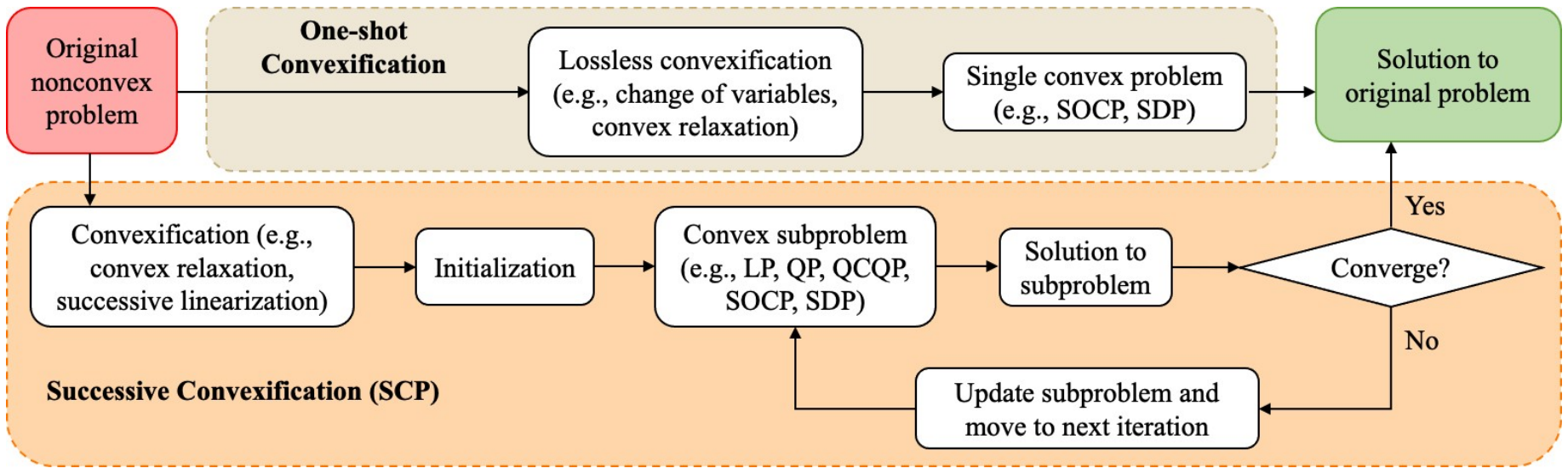

Figure 1: Schematic representation of one-shot lossless convexification and sequential convex programming (SCP) for nonconvex optimal G&C problems.

A general convex optimization problem takes the following form [10]:

*Problem 2*:

$$\begin{aligned} \text{minimize:}\quad & f_0(x) \\ \text{subject to:}\quad & f_i(x) \le 0,\ i = 1,2,\dots,m \\ & a_i^T x = b_i,\ i = 1,2,\dots,p \end{aligned}$$

where $f_0, f_1, \dots, f_m : \mathbf{R}^n \to \mathbf{R}$ are convex functions of $x \in \mathbf{R}^n$, which is a vector of the design variables and may represent the state and control sequences resulted from discretization. The convex optimization problem is essentially a special class of optimization problems with a convex objective function, convex inequality constraint functions, and affine equality constraint functions. The feasible set of a convex optimization problem is convex, and we minimize a convex objective function over this convex set. There are several sub-classes of convex optimization problems that are briefly summarized below.

### *2.2.1. Linear Programming*

When the objective and constraint functions are all affine, the problem becomes a linear programming (LP) problem that has the following form:

*Problem 3*:

$$\begin{aligned} &\text{minimize} \quad c^T x \\ &\text{subject to} \quad Ax = b \\ &\qquad\qquad\quad x \geq 0 \end{aligned}$$

where the parameters $c \in \mathbf{R}^n$, $A \in \mathbf{R}^{m\times n}$, and $b \in \mathbf{R}^m$ define the problem. The only inequalities are the componentwise nonnegativity constraints $x \geq 0$. The feasible set of the LP problem is a polyhedron, and the problem minimizes the affine function over this polyhedron.

*2.2.2. Convex Quadratic Programming*

If the objective in Problem 2 is a convex quadratic function and the constraint functions are all affine, the problem becomes a convex quadratic programming (QP) problem shown below:

*Problem 4*:

$$\begin{aligned} &\text{minimize} \quad \frac{1}{2} x^T P x + q^T x + r \\ &\text{subject to} \quad Ax = b \\ &\qquad\qquad\quad x \geq 0 \end{aligned}$$

where $P \in \mathbf{S}^n_+$ is symmetric positive semidefinite, $q \in \mathbf{R}^n$, $r \in \mathbf{R}$, $A \in \mathbf{R}^{m\times n}$, and $b \in \mathbf{R}^m$. In a convex QP problem, a convex quadratic function is minimized over a polyhedron.

*2.2.3. Convex Quadratically Constrained Quadratic Programming*

If the objective function and the inequality constraints in Problem 2 are all convex quadratic functions, the problem is a convex quadratically constrained quadratic programming (QCQP) problem as follows:

*Problem 5*:

$$\begin{aligned} &\text{minimize} \quad \frac{1}{2} x^T P_0 x + q_0^T x + r_0 \\ &\text{subject to} \quad \frac{1}{2} x^T P_i x + q_i^T x + r_i \leq 0, \; i = 1, 2, \ldots, m \\ &\qquad\qquad\quad Ax = b \end{aligned}$$

where $P_i \in \mathbf{S}^n_+$, $i = 0, 1, ..., m$ are symmetric positive semidefinite. The feasible region of a convex QCQP problem is the intersection of ellipsoids, and the problem is to minimize a convex quadratic function over this region.

*2.2.4. Second-Order Cone Programming*

LP, convex QP, and convex QCQP problems can be formulated as second-order cone programming (SOCP) problems, and SOCP is also a special case of convex optimization of the form:

*Problem 6*:

$$\begin{aligned} \text{minimize} \quad & f^T x \\ \text{subject to} \quad & \|A_i x + b_i\|_2 \le c_i^T x + d,\ i = 1,2,...,m \\ & Fx = g \end{aligned}$$

where $\|\bullet\|_2$ is the Euclidean norm. The problem parameters are $f \in \mathbf{R}^n$, $A_i \in \mathbf{R}^{n_i \times n}$, $b_i \in \mathbf{R}^{n_i}$, $c_i \in \mathbf{R}^n$, $d_i \in \mathbf{R}$, $F \in \mathbf{R}^{p \times n}$, and $g \in \mathbf{R}^p$. The second equation is called second-order cone constraint, and the feasible set of Problem 6 is the intersection of conic regions.

*2.2.5. Semidefinite Programming*

An SOCP can be formulated as a semidefinite programming (SDP) problem, which is a more general convex optimization problem of the following form:

*Problem 7*:

$$\begin{aligned} \text{minimize} \quad & c^T x \\ \text{subject to} \quad & x_1 F_1 + x_2 F_2 + \cdots + x_n F_n + G \preceq 0 \\ & Ax = b \end{aligned}$$

where $G, F_1, F_2, \ldots, F_n \in \mathbf{S}^k$ are symmetric $k \times k$ matrices, $A \in \mathbf{R}^{p \times n}$, and the inequality is a linear matrix inequality (LMI).

When $c_i = 0$, the SOCP problem is equivalent to a convex QCQP problem by squaring each of the constraints. Similarly, when $A_i = 0$, the SOCP problem reduces to an LP program. In addition, QCQP problems include QP problems as a special case by taking $P_i = 0$ in Problem 5, and QP problems include LP problems as a special case by taking $P = 0$ in Problem 4. As such, SOCP problems are more general than LP, convex QP, and convex QCQP problems.

A great deal of research in the convex-optimization-based G&C domain has been focused on the relaxation of the original OCP (Problem 1) into a convex optimization problem defined above and showing that an optimal solution to the relaxed convex problem is also an optimal solution to the original problem. Many nonconvex OCPs can be convexified by either restricting the original feasible set to a convex subset or enlarging the feasible set into a convex set containing the original feasible set [17]. However, proving the equivalence of the convexification process is not always possible, and guaranteeing that the convexification is lossless is difficult and highly problem-dependent. In fact, both approaches mentioned above can lead to some loss in the optimality or feasibility of the solution [17]. Lossless convexification emerged as a promising technique that allows obtaining the optimal solution to the original problem by solving a convexified, equivalent problem without removing the feasible region.

There are a few important theoretical contributions to the fundamental lossless convexification technique that relaxes particular types of OCPs into equivalent convex optimization problems, primarily focusing on OCPs with linear dynamics and nonconvex annular control constraints. For example, [58] considered a class of finite-time-horizon OCPs with continuous-time linear systems, convex cost, convex state constraints, but nonconvex control constraints. The control constraints were the only source of nonconvexity. A lossless convexification approach was proposed to relax the nonconvex control constraints, and the optimal solution to the relaxed problem was proved to be an optimal solution to the original nonconvex OCP. This lossless convexification approach was then extended to finite-time-horizon OCPs with continuous-time nonlinear dynamics and nonconvex control constraints [59]. Later, the convexification results were generalized to cases with additional linear or quadratic state constraints, where the convexification was still guaranteed to be lossless [60, 61, 62]. The proofs of these lossless convexifications were achieved using the maximum principle [63, 64], and the relaxed SOCP problems can be solved very efficiently. More recently, a class of mixed-integer nonconvex OCPs has been added to the list of problems that can be addressed by the lossless convexification approach, where the control input norms are restricted to be zero or lower- and upper-bounded [65]. Meanwhile, by removing some of the assumptions on system controllability and the gradient of the final point made in [58, 59, 62], more conditions on the validity of lossless convexification have been established for both free and fixed final time OCPs with nonconvex annular control constraints [66]. Furthermore, certain nonconvex OCPs with linear time-varying systems defined on disconnected control sets have shown to be potentially relaxed into convex problems using extreme point relaxations and normality approximations [67]. More theoretical advances in lossless convexification techniques are expected to emerge for relaxing more general OCPs into single convex optimization problems.

### *2.3. Sequential Convex Programming*

Given the fact that most of the G&C problems are not naturally in convex forms, series of transformation and relaxation techniques need to be employed to convert the original problem into a convex problem. For example, if highly nonlinear dynamic systems and nonconvex path constraints are incorporated into the problem, it may be difficult or impossible to formulate and solve a single convex optimization problem to find an optimal solution to the original problem. Instead, the nonconvex terms may be approximated by a successive process, in which the solutions of a sequence of convex subproblems are sought. This motivated the SCP method. As shown in Figure 1, the SCP method tackles a nonconvex OCP by repeatedly constructing and solving a convex subproblem in each iteration. The convex subproblem can be an LP, a convex QP, a convex QCQP, an SOCP, or less commonly an SDP problem. Each convex subproblem serves as an appropriate convexification of the original problem and is usually parameterized using the solution from the previous iteration. The process is repeated with an aim of making progress towards an optimal solution to the original problem [68]. More detailed discussion on SCP and its implementation can be found in [44].

It is worth noting that the SCP approach is similar to the popular sequential quadratic programming (SQP) algorithm, with the exception that more general set of convex optimization problems, such as QCQP and SOCP problems, are used as approximate subproblems during the iterative process. When the problem becomes highly nonconvex and cannot be handled as a convex problem, an SCP approach can be explored, in which the convex terms remain the same, but the nonconvex terms will be convexified through convex approximations of inequalities and affine approximations of equalities.

In the past 10 years, the SCP technique has advanced with a number of improvements for solving highly nonconvex G&C problems. The major differences in the SCP algorithms lie in what approximation approaches are used for convex relaxation, how the convex subproblems are formulated, how the intermediate solutions are used to parameterize the subproblem, what methods are used to measure the performance of the progress, and how theoretical guarantees can be enabled in terms of convergence and solution optimality. For example, [19] presented a successive SOCP-based convexification method for solving nonconvex OCPs with linear time-varying dynamics, and the nonconvexity arises from concave state inequality constraints and nonlinear terminal equality constraints. The concave inequality constraints were approximated by successive linearization, while the nonlinear equality constraints were handled by first-order expansions and compensated by second-order corrections. Guarantees were provided on the satisfaction of the original inequality constraints and the equivalence of the solutions to both the original problem and the converged successive solution. Successive convexification (SCvx) algorithms have also been developed to solve nonconvex OCPs in the presence of nonlinear dynamics and possible nonconvex state and control constraints [69, 70]. The nonconvex dynamics and constraints are generally convexified via successive linearization with respective to the solution of the subproblem solved in the previous iteration. However, an undesirable phenomenon has been observed in this process, i.e., an infeasible convex subproblem may be resulted in even if the original nonconvex problem is feasible. This phenomenon has been referred to as artificial infeasibility in the literature [69, 70]. With the aid of virtual control and trust region techniques, the artificial infeasibility issues can be mitigated, and global and super-linear convergence can be guaranteed under mild assumptions [69, 70].

Later, [71] generalized the earlier SCP-based methods in its guaranteed sequential trajectory optimization (GuSTO) framework for control-affine systems with drift, control and state constraints, and goal-set constraints under either fixed or free final time. The framework is guaranteed to converge to at least a stationary point. However, it was soon found that a general class of SCP-based methods are susceptible to an undesirable crawling phenomenon where slow convergence is observed when the algorithm is still far from a solution to the original nonconvex problem. This is usually the case when trust region and solution update rules with fixed iteration parameters are used to ensure feasibility and facilitate convergence. Potential remedies such as the use of hybrid algorithms are promising to mitigate this phenomenon [72]. Recently, the feasibility issues of the standard SCvx-based methods (i.e., the converged solution may not be feasible to the original nonconvex problem) has been addressed by incorporating the SCvx-based iteration in an augmented-Lagrangian-based framework [73]. In addition, the indirect method has been used to improve the performance of SCP-based methods for solving continuous-time OCPs with manifold-type constraints in [74], where convergence guarantees were established for control-affine dynamics. Techniques from the general numerical optimization field, such as line search and trust region, have also been employed to enhance the convergence of SCP in [27], where a line-search SCP algorithm and a trust-region SCP algorithm were developed. In the following years, it can be foreseen that the SCP-type methods will gain more popularity, and more theoretical analysis with more rigorous performance guarantees are expected to be reported for solving wider vehicular G&C problems.

## 3. Applications to Space Vehicles

G&C is a fundamental component for space vehicle systems and is crucial to the overall mission success. In this section, we will survey the applications of convex optimization for G&C of various space missions including powered descent guidance, small body landing, rendezvous and proximity operations, orbital transfer, space-craft reorientation, space robotics and manipulation, spacecraft formation flying, and orbital station keeping.

### *3.1. Powered Descent Guidance*

Rooted from the Apollo program [75], the powered descent guidance (PDG) technology has been one of the top priorities for a variety of manned and robotic space missions ranging from Moon landing to the exploration of Mars and other planets. More recently, PDG has gained additional interest and importance in the commercial space domain. The successful recovery and reuse of rocket boosters and stages by companies including SpaceX and Blue Origin has showed the great potential of PDG in reducing launch costs and improving mission responsiveness [76].

The PDG problem can be defined as generating an optimal trajectory that guides the vehicle from its initial or current state to a desired target state on the surface of the planet/moon with an expected accuracy of less than several hundred meters, which has been called precision landing or pinpoint landing (see Figure 2). However, solving a PDG problem is not easy. Many constraints, such as nonlinear vehicle dynamics, nonconvex constraints on the magnitude of the available thrust, and various state constraints, along with factors including mission uncertainties and environmental disturbances, give rise to a number of challenges when developing PDG systems and algorithms. Unlike lunar soft-landing where the problem has been well characterized and closed-form solutions have been obtained, the general three-dimensional (3-D) constrained PDG is much more difficult to solve. Also, since onboard computation of the flight trajectory and G&C command is expected, it is essential to exploit the structure of the problem and design algorithms with guaranteed convergence and real-time performance. These motivated the application of convex optimization for PDG problems. The representative works are summarized in Table 1 followed by a detailed literature review on this topic.

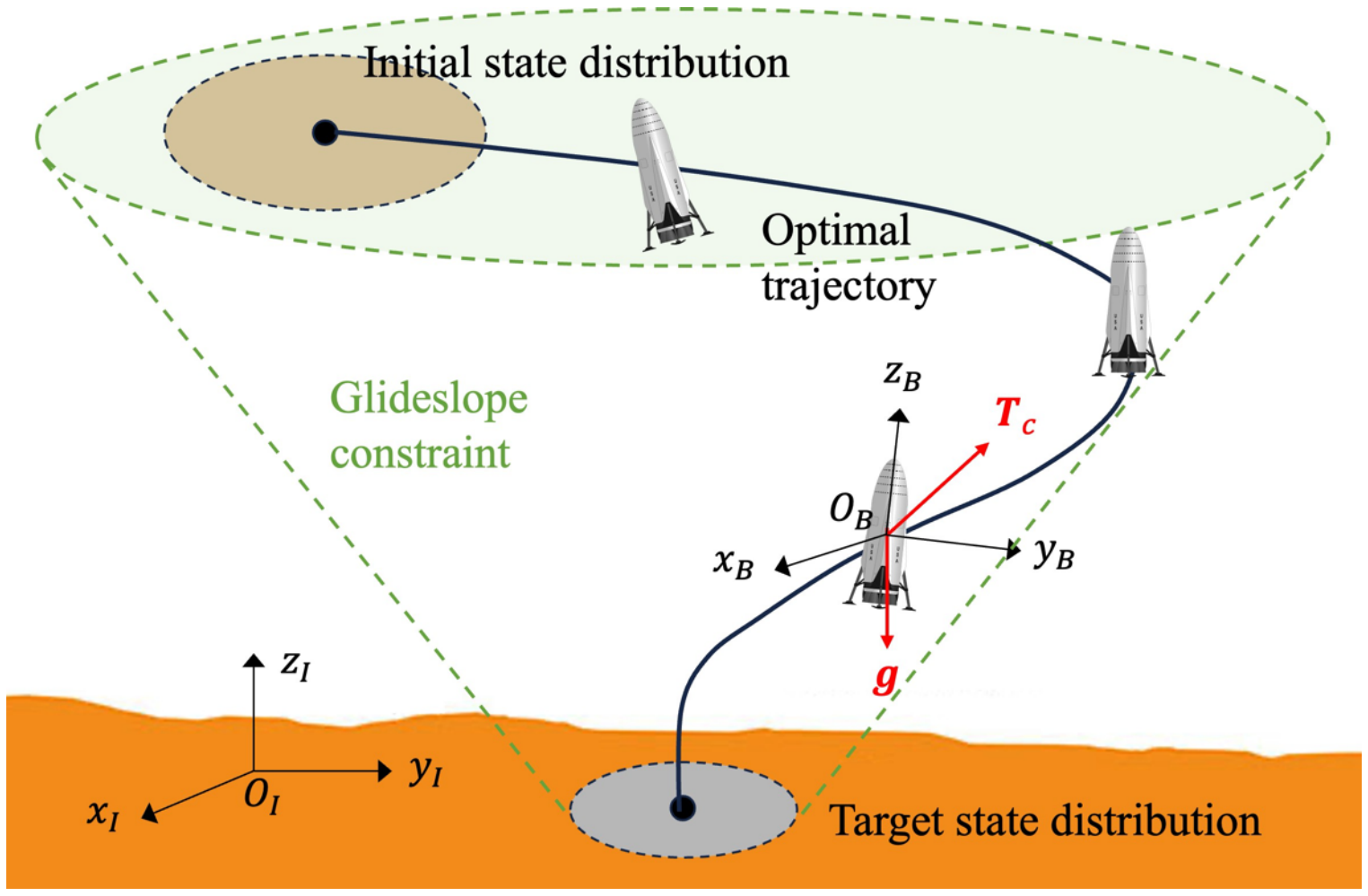

Figure 2: Schematic representation of powered descent guidance (PDG) for precision landing on the Moon, Mars, or other planets. The glideslope constraint ensures that the vehicle remains within a cone and stays at a safe distance from the ground until reaching the target.

Table 1: Summary of representative publications on convex optimization for powered descent guidance (PDG).

| Reference | PDG Problem | Approach | Slack Variables | Change of Variables | Convex Relaxation | Linear Approximation | Solver |
|---|---|---|---|---|---|---|---|
| [77] | 3-DoF fuel-optimal | SDP | ✓ | ✓ | ✓ | | SeDuMi |
| [15] | 3-DoF fuel-optimal | SOCP | ✓ | ✓ | ✓ | ✓ | SeDuMi |
| [16] | 3-DoF minimum-landing-error | SOCP | ✓ | ✓ | ✓ | ✓ | SeDuMi |
| [17] | 3-DoF with thrust pointing constraints | SOCP | ✓ | ✓ | ✓ | ✓ | SeDuMi |
| [78] | 3-DoF maximum-divert | SOCP | ✓ | ✓ | ✓ | ✓ | ECOS |
| [79, 80] | 3-DoF fuel optimal | SOCP + pseudospectral | ✓ | ✓ | ✓ | ✓ | ECOS, SDPT3 |
| [81] | 3-DoF multi-phase fuel-optimal | SOCP + pseudospectral | ✓ | | | ✓ | ECOS |
| [82] | 6-DoF fuel-optimal | SOCP + MPC | | ✓ | ✓ | ✓ | IPOPT |
| [83] | 6-DoF fuel-optimal | SOCP + SCP | ✓ | | ✓ | ✓ | SDPT3 |
| [84] | 3-DoF fuel-optimal with aerodynamic forces | SOCP + SCP | ✓ | ✓ | ✓ | ✓ | SDPT3 |
| [85] | 2-D fuel-optimal with thrust and aerodynamic controls | SOCP + SCP | | ✓ | ✓ | ✓ | MOSEK |
| [86, 87] | 6-DoF with state-triggered constraints | SOCP + SCP | ✓ | ✓ | ✓ | ✓ | SDPT3, ECOS |
| [88] | 3-DoF fuel-optimal with aerodynamic forces | SOCP + MPC | ✓ | | | ✓ | MOSEK |
| [89, 90] | 3-DoF fuel-optimal closed-loop | SOCP + covariance control | ✓ | ✓ | ✓ | ✓ | MOSEK |

The earliest work on convex-optimization-based PDG seems to be the one presented in [77], where an SDP-based solution was sought for fuel-optimal Mars pinpoint landing and showed better performance in terms of maneuvering than the polynomial guidance laws [91]. The results were then extended and published in [15], where the problem was relaxed into an SOCP problem, a subclass of SDP. In this pioneering work, the PDG problem was formulated as a fuel-optimal trajectory optimization problem subject to state and control constraints. The major nonconvex constraint lies in the nonzero lower bound of the thrust magnitude, i.e., $0 \leq \rho_1 \leq ||\boldsymbol{T}_c(t)|| \leq \rho_2$, which defines a nonconvex feasible region in the control space. To relax this nonconvex constraint, a slack variable was introduced, and it was proved that any optimal solution to the relaxed problem is also an optimal solution to the original problem. Therefore, the convexification was lossless. The relaxed OCP was finally discretized into an SOCP that can be solved onboard in real-time. This approach has been showed to have higher robustness than modified Apollo descent guidance algorithms and have stronger numerical stability and extensibility than constrained gradient-based indirect optimal control algorithms and analytic energy-optimal algorithms [92].

The convexification-based PDG algorithm developed in [15] has been extended to address problems considering more complex effects and constraints. For example, the algorithm was enhanced in [93] by including the rotation rate of Mars and extra state constraints into the problem formulation and introducing efficient ways to compute the optimal time-of-flight and detect the feasibility of the problem before solving it. Nonconvex attitude constraints due to thrust pointing have also been considered and expressed as $\hat{\boldsymbol{n}}^T \boldsymbol{T}_c(t) \geq ||\boldsymbol{T}_c(t)||\cos\theta$, which is convex for $\theta \in [0^\circ, 90^\circ]$ but nonconvex for $\theta \in (90^\circ, 180^\circ]$ [94]. An additional relaxation has been introduced to convexify the pointing constraint when $\theta \in (90^\circ, 180^\circ]$ such that the lossless convexification of the improved PDG algorithm remains for both the thrust bound and thrust pointing constraints [95]. The lossless convexification approaches presented in [15, 16, 95] were unified in [17], where both the minimum-fuel and minimum-landing-error PDG problems were solved under this unified optimization framework with thrust pointing constraints. Additionally, nonconvex obsta-cle avoidance constraints due to complex and hazardous terrains have also been incorporated in optimal PDG problems, which were converted into a sequence of convex subproblems via linearization, relaxation, and convexification techniques in [96, 97, 98]. More recently, the PDG approaches developed for Mars landings have been adapted to lunar soft landings. For example, a three-degrees-of-freedom (3-DoF) fuel-optimal PDG problem was established considering more state and control constraints including maximum tilt rate, maximum tilt acceleration, and maximum thrust ramp rate along with an inverse square gravity model and a minimum altitude constraint [99]. The lossless convexification framework from [15] was used to relax the problem into an SOCP problem.

PDG problems with different objectives have also been handled by convex optimization. For example, the distance to the prescribed target, $||\boldsymbol{r}(t_f)||^2$, was used as the objective to generate the minimum-landing-error trajectory via lossless convexification in [16]. Specifically, the algorithm determines the minimum-fuel trajectory to the target if a feasible trajectory exists; if no feasible trajectory to the target exists, however, it calculates the trajectory that minimizes the landing error. SOCP was used in both situations. In addition, the results in [77] and [15] have been extended to address the maximum-divert planetary landing problem with linear and quadratic state constraints such as the velocity constraints imposed to keep bounded aerodynamic forces and ensure the structural integrity of the vehicle [78]. Due to this change, new theoretical results have been derived to ensure lossless convexification for flights where these state constraints are active. To improve the landing accuracy, a navigation-optimal PDG problem has been solved via successive linearization of the covariance matrix elements of the Extended Kalman Filter (EKF) of the vehicle's navigation algorithm by minimizing the trace value of the covariance matrix at landing [100]. More recently, the method has been extended to solve the minimum-landing-error PDG problem under stochastic navigation errors that were modeled as chance constraints [101]. The design of PDG algorithms under uncertainties and disturbances will be discussed shortly below.

Discretization is crucial for the convex-optimization-based direct method that converts the original continuous-time problem into a discrete parameter optimiza-tion problem. A proper discretization can decrease the time required to find a solution of acceptable accuracy while satisfying the real-time computational requirement. Aiming to improve the accuracy of the convex approach without excessively worsening its real-time performance, pseudospectral methods and convex optimization have been combined to solve OCPs such as PDG. For example, the flipped Radau pseudospectral method and the Lobatto pseudospectral method have been combined with convex optimization, leading to a flipped Radau pseudospectral convex method and a Lobatto pseudospectral convex method, respectively [79]. These two methods have been applied to solve the fuel-optimal Mars PDG problem and compared with the standard convex method in [79], where the reported simulation results showed that the pseudospectral convex methods are capable of producing more accurate results than the standard transcription methods such as finite differences and the trapezoidal rule. When the number of nodes grows, however, the solution time with be larger than the standard methods due to the loss of sparsity in the discretized optimization problem [102]. This issue can be mitigated by generalizing the standard pseudospectral-convex method in the frame of the broader family of *hp* schemes and developing a hybrid framework consisting of *hp*-pseudospectral methods and convex optimization to significantly reduce the computational time [80]. In combination with SCP, the performance of pseudospectral methods have been further compared with other discretization methods including zero-order hold (i.e., keep the control input constant between sampling times), first-order hold (i.e., define the control input as a linear function between sampling times), and the classical fourth-order Runge-Kutta method in solving a fuel-optimal PDG problem [102]. The results suggested that pseudospectral methods are capable of producing more consistent trajectories and less sensitive to the discretization resolution than other discretization methods. In addition, pseudospectral method and convex optimization have been combined in solving 3-DoF multi-phase fuel-optimal PDG problems [81], where the optimal phase division was determined by the indirect method. Recently, the Radau pseudospectral method has been combined with convex optimization for solving a 2-D vertical landing problem of a starship-like vehicle with large attitude flip, and the effectiveness of the method has been verified via hardware-in-the-loop experiments [103].

Comparing to the 3-DoF PDG problems, the 6-DoF PDG scenarios are more challenging to address due to the thrust-magnitude lower bound, the mass depletion dynamics, plus the additional highly nonlinear attitude dynamics. A single convex optimization problem cannot be formulated in this case to find the solution to the original problem; instead, a successive convexification can be employed to handle these nonconvexities. For example, a fixed-final-time 6-DoF fuel-optimal Mars PDG problem has been solved in [83], where the original problem was transformed into a sequence of SOCP problems. To facilitate convergence, quadratic trust regions were introduced to keep the solution bounded and a relaxation term was added to the dynamics to ensure feasibility throughout the convergence process. Later, this successive convexification framework was extended to solve the minimum-time 6-DoF PDG problem where the time of flight is free to be minimized subject to similar constraints [104]. Following an unconventional means of representing the orientation and position of the lander spacecraft, dual quaternions have been used to simultaneously represent the rotational and translational motion dynamics [105, 82]. One particular feature about this dual-quaternion parameterization method is that the equations of motion can be expressed in a form similar to the standard quaternion kinematic and dynamic equations, and some constraints such as line of sight can be expressed in convex forms over a given set of dual quaternions [106]. By leveraging this attractive feature, the fuel-optimal PDG problem has been solved within the framework of piece-wise affine MPC, where the resulting nonconvex constraints (e.g., line-of-sight constraints and glideslope constraints) were converted into computationally tractable convex constraints for onboard computation [105, 82]. The benefits of combining convex optimization with MPC to solve PDG problems will be discussed in more detail below.

It has been observed that aerodynamic forces were generally neglected in the development of SOCP-based PDG methods reviewed above. Incorporating aerodynamic forces in the problem formulation would potentially result in more practical solutions but add significant complexity to the problem such that the previously reviewed SOCP-based algorithms may fall short in finding accurate optimal solutions.

This type of problems falls into the scope of atmospheric flight missions that will be discussed in Section 4; however, publications on this topic are reviewed here as a natural extension of the PDG problem. [84] seems to be the first work that incorporated aerodynamic forces in the development of convex optimization approaches for PDG problems. The nonlinearities introduced by aerodynamic drag, mass-depletion dynamics, and free time-of-flight cause critical challenges for real-time applications. In addition to utilizing the lossless convexification method to address the minimum-thrust constraint, successive convexification relying on the use of linearization, trust regions, and relaxations has been employed to eliminate the remaining nonconvexities, leading to a sequence of iteratively solved SOCP problems. The method has been showed to converge in a small number of iterations and robust to a wide range of time-of-flight guesses [84]. Following a similar lossless convexification approach, a 3-DoF fuel-optimal rocket landing problem considering aerodynamic drag and Earth rotation has been solved in [107], where an improved successive convexification method was developed by leveraging the efficient pseudospectral method and a dynamic trust-region updating strategy. Other than using only the thrust magnitude and thrust direction as the controls, both aerodynamic forces and engine thrust can be used as control inputs for PDG. For example, [85] addressed a 2-D PDG problem for a reusable rocket returning back to Earth by coordinating the thrust and the aerodynamic forces to achieve fuel-optimal landing. The nonconvex constraints on the aerodynamic forces and the thrust were converted into convex forms using relaxation and linear approximation. As a result, the original problem was transformed into and solved as a sequence of SOCP problems. Following the lines in [84], the successive-convexification-based algorithm has been tailored to PDG of reusable launchers over extended flight envelopes, and the algorithm has been implemented and verified in a closed-loop manner on a fuel-optimal PDG problem with aerodynamic and thrust forces in [108]. A 6-DoF vehicle model was formulated in this work; however, only the 3-DoF equations of motion were considered in the optimal control problem formulation and the development of the successive SOCP algorithms under a perfect attitude control assumption. More recently, the successive convexification approach has been applied to PDG of parafoil for precision landing on Titan, Saturn's largest moon with a dense atmosphere [109]. By using only aerodynamic forces for control, a 6-DoF minimum-control-effort parafoil PDG problem was formulated and solved as a sequence of SOCP problems via successive convexification with the aid of flexible trust regions and virtual controls.

In addtion, optimal PDG considering state-triggered constraints has been an active research area in recent years due to the fact that PDG missions usually involve multiple flight phases such as braking, approach, and final descent, and some constraints need to be enforced only when certain criteria are satisfied. Such constraints include state-based keep-out zone constraints for collision avoidance and distance-based line-of-sight pointing constraints. Fortunately, these constraints can be incorporated in the problem formulation while maintaining the continuity of the optimization framework without introducing binary/integer variables and resorting to time-consuming heuristics or mixed-integer programming algorithms. One of the first results on state-triggered PDG was reported in [106], where a line-of-sight pointing constraint was enforced based on the distance from the landing site. Following the dual-quaternion parameterization in [105, 82], a 6-DoF nonconvex fuel-optimal PDG problem was formulated and solved using a successive convexification procedure. The optimal guidance trajectories were generated with the line-of-sight constraint explicitly enforced when the trigger condition is satisfied [86]. A further step has been taken in [110] and [87] by generalizing the state-triggered constraint formulation into compound state-triggered constraints defined by vector-valued trigger and constraint functions using Boolean logic operations. A continuous formulation for these compound state-triggered constraints was established and handled by successive convexification in [110], which applied the approach to a 6-DoF minimum-time PDG problem with an ellipsoidal aerodynamic model and a free-ignition-time modification. A velocity-triggered angle-of-attack constraint and a collision-avoidance constraint were considered as two examples of such compound state-triggered constraints. Later, the convergence of the successive convexification approach was enhanced by virtue of virtual control and trust region modifications in [87], where a free-final-time fuel-optimal 6-DoF PDG problem was solved with velocity-triggered angle-of-attack and range-triggered line-of-sight constraints.

Furthermore, earlier publications on convex-optimization-based PDG generally assumed a deterministic OCP formulation, and the solution usually resulted in a nominal trajectory that connects a single pair of initial and terminal states without considering uncertainties or external disturbances. However, many factors, such as modeling uncertainties, localization errors, and environmental disturbances may cause substantial deviations from the nominal optimal trajectory during actual flight. One solution to this problem is to continuously update and follow the nominal trajectory online through closed-loop control [111]. MPC is such an approach to achieving this purpose. By recursively solving constrained optimization problems online with the repeatedly updated system states, MPC is robust to uncertainties and disturbances during the flight and has found wide applications for G&C of vehicular systems [54, 55]. Convex optimization and efficient discretization methods such as pseudospectral collocation can be combined and implemented within the MPC framework to develop receding-horizon PDG schemes with a certain degree of robustness. For example, a 3-DoF fuel-optimal PDG problem has been formulated in [88] considering both aerodynamic and thrust controls, and a pseudospectral-based successive convexification algorithm was used to solve the problem under an MPC framework to rapidly compute the optimal trajectory in each MPC circle with a high trajectory update frequency. To achieve more precision prediction and better control, a 3-D special Euclidean group, SE(3), has been used to establish the 6-DoF vehicle dynamics, which has been discretized using a Lie group variation integrator from geometric mechanics, leading to accurate yet robust algorithms for PDG in combination with convex relaxation and nonlinear MPC [112].

An alternative approach to addressing uncertainties and disturbances is to explicitly account for the stochasticity in PDG formulations, leading to better performance than simply employing deterministic closed-loop control [113]. To this end, a stochastic extension to the deterministic convex-optimization-based PDG has been studied for more general PDG scenarios where the vehicle is steered from some initial state distributions to some target state distributions with Brownian motion process noise acting on the system as external disturbance forces [114]. Then, the gen-eration of the nominal trajectory is coupled with the design of the closed-loop control law through a stochastic PDG problem formulation that can be solved through optimal mean control and optimal covariance control [115, 116, 117]. By modeling the PDG as a stochastic process, this approach is capable of obtaining less conservative feed-forward optimal thrust command to allow for sufficient feedback authority [89]. More recently, by characterizing vehicle's mass as a random variable with variance, a stochastic 3-DoF fuel-optimal PDG problem has been formulated as a covariance control problem by explicitly including uncertainties and disturbances in the formulation [90]. Through a convenient change of variable and successive convexification, the covariance control problem has been cast as a sequence of deterministic convex optimization programs, from which the optimal nominal trajectory and the feedback control policy can be obtained simultaneously. Moreover, with the same goal of developing fast yet robust approaches to PDG, the polynomial chaos theory has been combined with convex optimization in solving stochastic optimal PDG problems in [118], where the polynomial chaos method was utilized for dynamic uncertainty propagation to calculate the mean and variance of the states, constraints, and performance index, thus transforming the original stochastic problem into a deterministic version in a higher-dimensional space. The transformed deterministic was then solved by successive convexification.

Physical experimentation and flight tests are important steps towards future onboard applications of the methods. The lossless convexification methodology and the associated SOCP-based algorithm have been integrated as the Guidance for Fuel Optimal Large Diverts (G-FOLD) tool [119]. Considering vehicle dynamics and relevant mission constraints, flight tests on JPL's Autonomous Descent and Ascent Powered-flight Testbed (ADAPT) have been performed using the Masten Space Systems Xombie vertical-takeoff vertical-landing suborbital rocket to demonstrate the off-line G-FOLD generated trajectories during the summer of 2012 [120], test the divert trajectories calculated by G-FOLD onboard in 2013 [121], and test the integration of the Lander Vision System and G-FOLD in two successful free flight demonstrations on the Xombie vehicle in December 2014 [122]. These results demonstrated that G-FOLD is capable of planning optimal trajectories respecting all the constraints of the rocket-powered vehicle. It is worth mentioning that the flight software successfully validated in these flight tests was based on the customized code for generation of optimal landing trajectories onboard in real-time [123]. By making use of sparsity, explicit code generation, and exact memory allocation, the IPM has been tailored for SOCP problems, producing customized ANSI-C code for embedded real-time applications. This was claimed to be the first real-time embedded convex optimization algorithm used to control large vehicles such as the ADAPT guided rocket [123]. This customized solver has been proved to be capable of providing accurate results rapidly enough for real-time applications through comparisons with generic solvers such as SDPT3, SeDuMi, and ECOS. It has been shown that per time of flight, infeasibility or an optimal trajectory can be calculated in approximately 0.7 s on a state-of-the-art radiation-hardened flight processor or in approximately 2.5 s when running in the background on a flight processor [124]. Later, an SCP algorithm has been refined to be compatible with common flight code requirements while maximizing its computational performance [125]. The reader is referred to [126] for a system design of ADAPT, a detailed description of the customized algorithmic components of the flight software implemented on ADAPT, the results of all the three years of flight tests as well as more in-depth analyses of implementation issues. Finally, as one of the core Precision Landing and Hazard Avoidance (PL&HA) technologies and capabilities, the 6-DoF dual-quaternion convex-optimization-based PDG algorithm [106, 102, 127] has been supported by NASA's Safe and Precise Landing – Integrated Capabilities Evolution (SPLICE) project for future robotic science and human exploration missions to the Moon, Mars, and other solar system bodies [128, 129]. Through rapid prototyping and coding, the SPLICE team has deployed a flight code version of this PDG algorithm on the descent and landing computer for test benchmarking [130].

### *3.2. Small Body Landing*

With the rapid development of space technologies, there has been an growing interest in exploration of small celestial bodies such as asteroids and comets. Over the years, small body exploration missions have gradually transformed from fly-by and orbiting to proximity operations (e.g., descent, landing, hopping), impact, and sample return missions [51]. Compared to large planetary bodies, however, small body missions face a number of unique challenges. For example, the gravity field of a small body is generally weak and difficult to accurately model due to its small size, irregular shape, non-uniform mass distribution, complex rotational state, and limited ground observation [131]. As a result, perturbations such as solar radiation and gravity from other celestial bodies may dominate the forces acting on the spacecraft and lead to additional uncertainties and disturbances in the dynamic model, which will degrade the mission performance [50]. Furthermore, unexpected obstacles due to the complicated terrain of small bodies and unexpected disturbances such as outgassing activities of comet-like celestial bodies will cause additional difficulties to proximity operations and surface exploratory missions [51].

Descent and landing on the target is one of the most critical phases for successful small body missions. Among the enabling technologies for safe landing, advanced G&C techniques are of paramount significance in determining mission success. A small body landing (SBL) G&C scheme determines appropriate control actions and the corresponding state trajectory that leads the spacecraft from its current state to the desired target landing state within a reasonable time at acceptable fuel consumption with minimum landing error subject to the highly nonlinear dynamics and the state and control constraints during the maneuver (see Figure 3) [132]. The produced control and state profiles can be used either as reference trajectories for an outer-loop controller to track or as the basis of MPC implementation [133]. However, due to the limited knowledge of the target, complex terrain environment, significant model uncertainties, external disturbances, and long-distance communication delays, the G&C system is expected to have some degree of autonomous capabilities. Specifically, the G&C algorithms are required to be computationally efficient and robust to model uncertainties and disturbances while incorporating state and control constraints for safe and rapid onboard decision-making. The current onboard computational power and algorithmic advances such as convex optimization make such G&C approaches possible. The literature on this topic is reviewed in detail below, and the representative publications are summarized in Table 2.

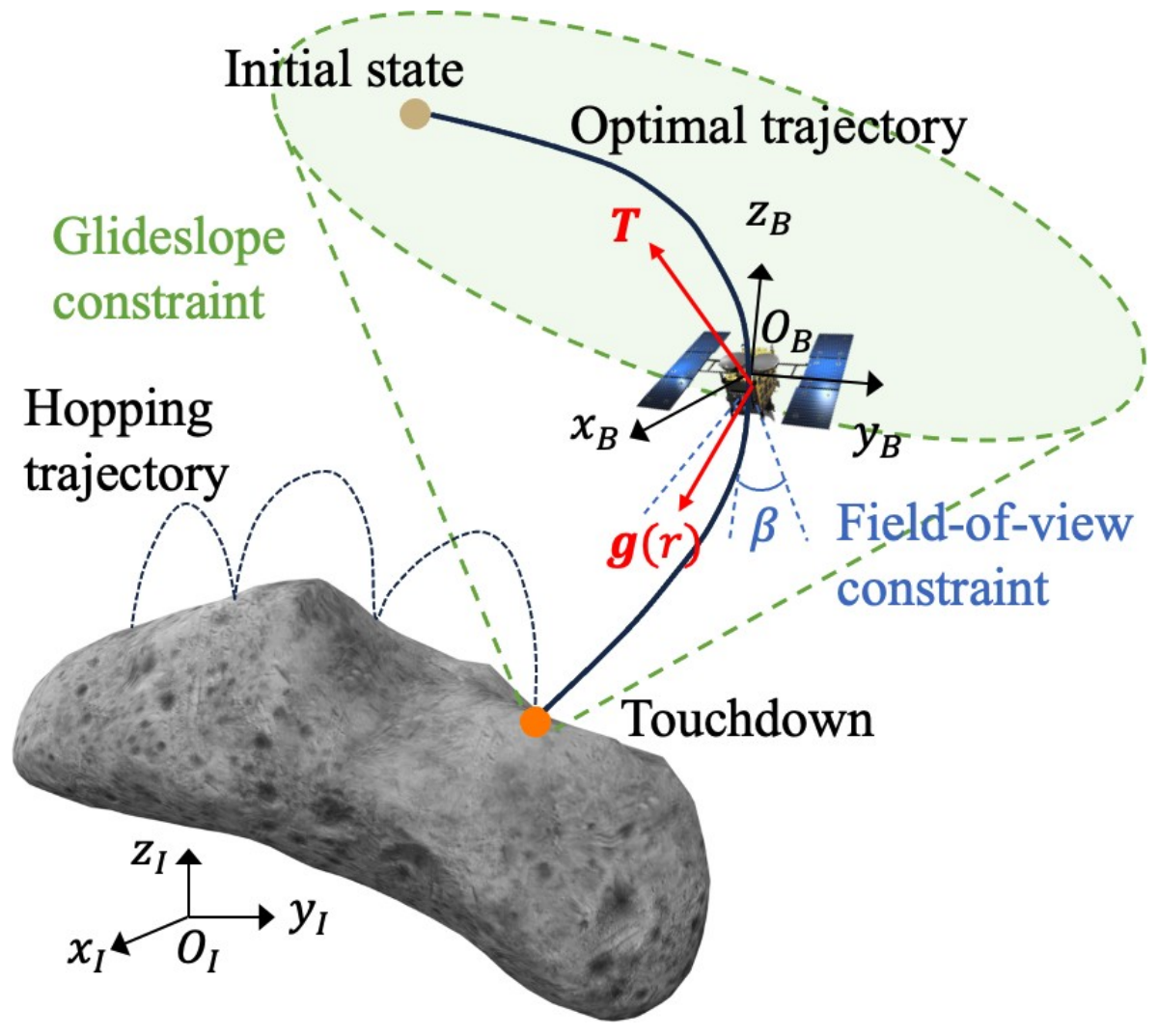

Figure 3: Schematic representation of small body landing (SBL).

Table 2: Summary of representative publications on convex optimization for small body landing (SBL).

| Reference | SBL Problem | Approach | Slack Variables | Change of Variables | Convex Relaxation | Linear Approximation | Solver |
|---|---|---|---|---|---|---|---|
| [134, 135] | 3-DoF fuel-optimal | SOCP + SCP | ✓ | ✓ | ✓ | ✓ | CVX |
| [136] | 3-DoF time-optimal | SOCP + SCP | | | ✓ | ✓ | SDPT3 |
| [137, 138] | 3-DoF fuel-optimal | SOCP + SCP | ✓ | ✓ | ✓ | ✓ | MOSEK |
| [139] | 6-DoF time-optimal with fuel consumption penalty | SOCP + SCP | | | | ✓ | SDPT3 |
| [140] | 6-DoF fuel-optimal | SOCP + SCP | ✓ | | | ✓ | SDPT3 |
| [141] | 3-DoF minimum-landing-error | SOCP + SCP | | ✓ | | ✓ | CVX |
| [132, 133] | 3-DoF minimum-control | SOCP + MPC | | | | ✓ | SeDuMi, SDPT3 |
| [142, 143] | 3-DoF multi-objective | QP + MPC | | | | ✓ | CVX |
| [144] | 3-DoF fuel-optimal | SOCP + MPC | ✓ | ✓ | ✓ | ✓ | Gurobi |

In light of the great success of convex relaxation and successive convexification techniques in planetary PDG, the applicability of these techniques to small body powered descent and landing problems has been investigated. The earliest work observed is [134] where a 3-DoF fixed-final-time fuel-optimal asteroid powered descent problem was solved to generate the optimal thrust profile for soft landing at the target location considering a gravity model of higher fidelity than Newtonian and various mission and operational constraints. Motivated by the techniques for planetary PDG [15, 16], a slack variable was introduced to relax the nonconvex control constraint (i.e., non-zero lower bound of the magnitude of the thrust vector) into a convex form consisting of $||\boldsymbol{T}(t)|| \leq T_m(t)$ and $T_{min} \leq T_m(t) \leq T_{max}$. The equivalence of this relaxation technique was theoretically established without including the glide slope constraint in the analysis. Then, the problem was cast into an SOCP problem via a change of variables and Taylor series expansions, and an optimal solution to the original problem was obtained through SCP with a successive approximation of the nonlinear gravitational acceleration term in the dynamics. The approach developed in [134] for triaxial ellipsoidal asteroids has been extended for landing on irregularly shaped asteroids in [135], where a higher-fidelity gravity model that balances accuracy and computational complexity was used, and a single-variable outer optimization loop was employed to find the optimal flight time that yields the overall best fuel economy.

Compared to the fuel-optimal SBL studied in [134] with a fixed flight time, the time-optimal SBL adds an additional nonconvex factor, i.e., the free final time, to the problem formulation. Instead of transforming the problem into an SOCP using the convexification techniques in [134], the time-optimal problem has been solved through connecting with other related problems. As mentioned before, the minimum-landing-error problem has been used to find a solution when no feasible solution exists for fuel-optimal Mars landing [16]; in contrast, the minimum-landing-error problem has been connected with the time-optimal SBL in [136]. First, a reduced minimum-landing-error SBL problem was solved using convex optimization based on the observation that the thrust remains at its maximum magnitude during the entire flight when the flight time is less than or equal to the minimum flight time.

Then, the minimum flight time was sought through a combination of extrapolation and bisection methods, instead of the line search method, to speed up the search [136]. It is worth noting that the lower bound constraint on the thrust magnitude was ignored in [136] due to the fact that the thrust must stay at its maximum magnitude for the solution to the time-optimal problem, which eliminate the need of convexification techniques studied in [134].

In addition to the descent and landing missions, hopping trajectories have also been studied via convex optimization for surface exploration on small bodies [137]. Trajectories for both single-hopping and multi-hopping scenarios can be generated for a hopper or a surface explorer to reach specific targets to perform exploration tasks. Through the similar convexification and relaxation techniques invested in [134, 15], the feasibility of using SOCP and SCP to generate fuel-optimal hopping trajectory has been demonstrated considering conic constraints on both the endpoints of the hopping trajectory [137]. Shortly after, this convex optimization approach was used to facilitate active trajectory control and more intelligent SBL strategies for both landing and hopping explorations on small bodies [138], where a new discretization method based on the explicit fourth-order Runge-Kutta rule was developed to improve the solution accuracy while maintaining real-time performance. Moreover, the hopping sequence can be optimized to aid the SCP-based long-distance hopping transfer on the asteroid surface [145].

In recent years, more results have been reported on using convex optimization for solving SBL problems with growing complexity. For example, a 6-DoF time-optimal SBL problem augmented by a fuel consumption penalty with two-phase free final time has been solved via successive convexification considering a glide-slope constraint for collision avoidance and an attitude constraint for field-of-view of the landing camera plus the constraints on the thrust, torque, and the mass of the vehicle [139]. Different from the combination of extrapolation and bisection methods used in [136] and the outer optimization loop employed in [135], the ideas of normalizing the flight time into the range of [0, 1] with time dilation coefficients from [104] was used in [139] to directly generate the optimal flight time along with the trajectory. Techniques such as virtual control and trust regions have been used to address the artificial infeasibility problems to promote the convergence of the successive solution process. Later, this approach was further developed to solve a 6-DoF fuel-optimal SBL problem in [140] in combination with the relaxation techniques in [15, 135] due to a slight modification of the problem formulation, i.e., the addition of a nonconvex lower bound constraint of the thrust magnitude. More recently, convex optimization has been used to address SBL of high area/mass ratio landers controlled by solar radiation pressure (SRP) [141]. Different from the conventional thrust-driven landers, SRP-propelled SBL suffers from new challenges due to the fewer degrees-of-freedom SRP control (i.e., controlled by only two angles). Through a successive convexification approach, the problem was converted into a sequence of SOCP problems, and a trust region constraint plus a modified objective function were introduced to improve the convergence of the SCP process [141].

To handle model uncertainties (e.g., errors in gravity model of the small body) and exogenous disturbances, convex optimization has been combined with MPC in developing robust SBL G&C algorithms. For example, linearized models of gravity have been used to formulate a linear time-varying model of SBL dynamics, based on which an SOCP-based 3-DoF minimum-control SBL guidance scheme with state and control constraints has been developed to facilitate real-time generation of open-loop pseudo way-point trajectories that can be updated and tracked in a robust MPC manner [132]. The approach was then augmented in a two-mode scheme in [133], where a standard mode guides the vehicle toward the desired target state in a receding-horizon manner, and a safety mode can be activated due to invalid expected state constraints and errors in state determination to maintain the vehicle at a safe state from the surface, providing some form of state-constraint robustness and risk mitigation. Furthermore, a convex MPC method has been developed for asteroid landing based on the linearized model in [142], where the landing mission was split into a circumnavigation and a landing phase, and a convex QP problem was solved at each time subject to linear equality constraints and affine hyperplane inequality constraints for collision avoidance. However, replacing the concave collision-avoidance constraint with affine approximations may lead to conservativeness. Also, the rate of hyperplane rotation is a design parameter that may be problem-dependent and difficult to choose. To address these issues, an optimal hyperplane method has been developed by solving a separate convex problem to free the vehicle from colliding with the surface of the asteroid [143]. In combination with a projection theorem argument, the vehicle can be guaranteed to converge to the desired target state. More recently, by integrating the convexification techniques in [134], change of variables in [15], and the two-phase approach in [142], a 3-DoF fuel-optimal SBL problem has been solved in [144] under an MPC framework to cope with unmodeled dynamics and disturbances during SBL maneuvers. In addition to the predictive controllers, other control strategies such as input observers and extended command governors have also been used to compensate gravity model errors and enforce state and control constraints for autonomous SBL [146]. The extended command governor has been shown to provide better fuel economy, while the MPC methods offer superior constraint handling and disturbance rejection performance at the expense of increased difficulty in tuning [142].

### *3.3. Rendezvous and Proximity Operations*

Space missions involving two or more vehicles have received increasing attention in recent years but also raised critical G&C challenges, especially for rendezvous and proximity operations (RPOs) [147]. Experiences have suggested that autonomous RPOs would greatly benefit from highly efficient G&C techniques with additional safeguards to protect the vehicles from potential mission failures [148]. In general, RPOs refer to the maneuvers of a chaser vehicle to approach an in-orbit target vehicle for missions such as flying around or docking with the target vehicle. An RPO mission may involve multiple phases, and this subsection focuses on the phases where the chaser has already arrived to the vicinity of the target. To enable autonomous RPOs, it is of imperative importance to develop mathematically rigorous G&C algorithms that can plan the mission trajectory and generate G&C commands reliably in real-time onboard the vehicle with minimum crew intervention or ground support [18]. A successful G&C method should be able to produce an appropriate trajectory and the associated control actions (usually thrust) that will lead the chaser to transit from its initial state to the target condition within a certain period of time while consuming as little fuel as possible and satisfying all the mission constraints such as approach corridor, sensor field-of-view, collision avoidance, plume impingement (see Figure 4).

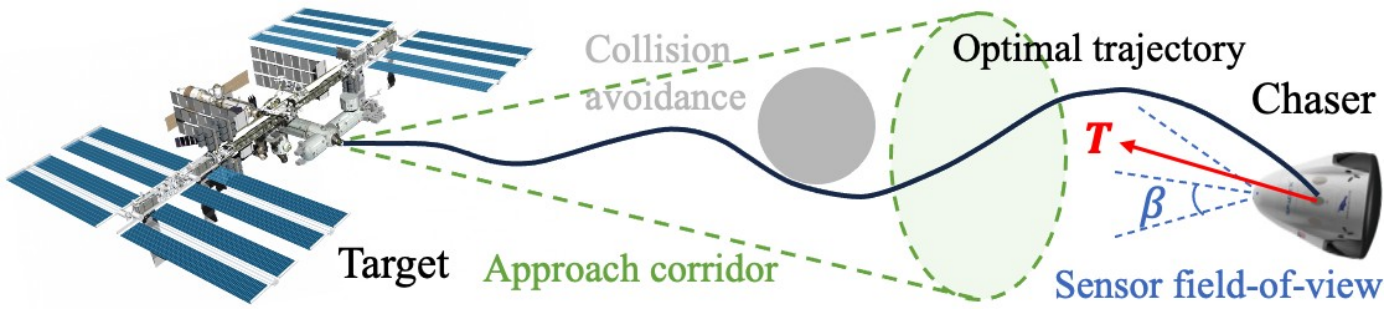

Figure 4: Schematic representation of rendezvous and proximity operations (RPOs).

Numerous publications on using optimal control and numerical optimization for generating RPO trajectories have been found in the literature for a wide range of RPO scenarios. Earlier relevant approaches, such as mixed-integer linear programming (MILP) [149, 150], MILP-based MPC [151, 152, 153], and linear-quadratic MPC [154, 155, 156, 157, 158, 159, 160], have proved to capable of addressing nonconvex RPO problems with various state and control constraints (e.g., collision avoidance, plume avoidance, line of sight) based on linearized system dynamics. Also, robust controllers have been designed via a Lyapunov approach for RPOs where the controller design problem can be cast into a convex optimization problem subject to liner matrix inequalities (LMIs) considering linear relative equations, parameter uncertainties, external perturbations, and control constraints [161, 162, 163, 164]. In addition, polynomial optimization has been combined with SDP and LMI-based convex relaxation in solving RPO problems [165, 166, 167]. Nearly all these works are based on linearized dynamics such as the Clohessy–Wiltshire (CW) equations that are easy to use to analyze and visualize relative RPOs but may be conservative in mission designs due to the restrictive assumptions of near-circular orbits of the target [18]. Therefore, there is a need for novel computational G&C approaches to RPOs on arbitrary orbits under nonlinear gravity models and potential perturbations such as $J_2$ harmonics, aerodynamic drag, and solar radiation pressure.

In the review below, we focus on the papers published over the past 10 years reporting more systematic development of convex optimization methods for RPOs via techniques such as lossless and successive convexification (see Table 3), not simply relying on linearization or the assumption that the target is in a circular orbit.

Table 3: Summary of representative publications on convex optimization for rendezvous and proximity operations (RPOs).

| Reference | RPO Problem | Approach | Slack Variables | Change of Variables | Convex Relaxation | Linear Approximation | Solver |
|---|---|---|---|---|---|---|---|
| [18, 168] | 3-DoF fuel-optimal | SOCP + SCP | ✓ | ✓ | ✓ | ✓ | MOSEK |
| [169] | 3-DoF multi-objective | MISOCP + MPC | ✓ | ✓ | ✓ | ✓ | MOSEK |
| [170] | 3-DoF minimum-time using differential drag | LP | | | ✓ | | CVXGEN[171] |
| [172] | 3-DoF fuel-optimal | SOCP + SCP | | ✓ | ✓ | ✓ | Gurobi |
| [173] | 3-DoF fuel-optimal and minimum-time | SOCP + SCP + mesh refinement | ✓ | ✓ | ✓ | ✓ | ECOS |
| [174] | 3-DoF fuel-optimal | SOCP + SCP + pseudospectral | ✓ | ✓ | ✓ | ✓ | MOSEK |
| [175, 176] | 3-DoF fuel-optimal | SDP | ✓ | ✓ | | | SDPT3 |
| [177] | 3-DoF fuel-optimal | SOCP + SCP | ✓ | ✓ | ✓ | ✓ | - |
| [178] | 3-DoF fuel-optimal | SOCP + SCP | | ✓ | ✓ | ✓ | MOSEK |
| [179] | 6-DoF fuel-optimal | SOCP + SCP | | ✓ | | ✓ | SDPT3 |
| [180, 181] | 6-DoF fuel-optimal with state-triggered constraints | SOCP + SCP | ✓ | | | ✓ | MOSEK, SDPT3 |
| [182] | 6-DoF fuel-optimal closed-loop | Covariance control + SDP + SCP | ✓ | ✓ | ✓ | ✓ | MOSEK |

Inspired by the success of lossless convexification and SOCP-based methods in PDG [15, 16], similar techniques have been pursued to solve RPO problems. As one of the earlier publications in this area, [18] posed the RPO problem as a 3-DoF fuel-optimal fixed-final-time OCP subject to an inverse-square gravity model, intrinsic nonlinear thrust terms in the equations of motion (as in [15, 16]), and trajectory constraints on approach corridor, thrust plume direction, terminal conditions, and possible intermediate way-points. Convex relaxation techniques similar to those in [15, 16] have been used to transform the original problem into a relaxed one. Specifically, a slack variable $\eta$ was introduced to place $||\boldsymbol{T}(t)||$ in the problem formulation. As a result, the original thrust control constraint, $||\boldsymbol{T}(t)|| \leq T_{max}$, became two relaxed ones, $||\boldsymbol{T}(t)|| \leq \eta(t)$ and $0 \leq \eta(t) \leq T_{max}$, and the nonconvex thrust direction constraint was replaced by a convex inequality constraint. The equivalence of the solutions to the RPO problems before and after the relaxation has been established [18]. To remedy the nonlinearity in the dynamic equations, the thrust acceleration was designated as part of the control through a change of variables, and the nonlinear gravity model was circumvented via a successive linearization method. Finally, the solution to the relaxed (and the original) nonconvex RPO problem was successively approached by the solutions to a sequence of SOCP problems with linear, time-varying dynamics. It is worth mentioning that this methodology is general and allows incorporation of more constraints including nonconvex keep-out zones and more complicated factors such as $J_2$ terms and aerodynamic drag, which are important for RPOs in low Earth orbits [168].

The methods in [18, 168] have inspired the development of similar approaches to collision avoidance maneuver optimization [183] and other RPO scenarios [184]. For example, the techniques in [18, 168] have been extended to solve RPOs with obstacle-avoidance constraints by reformulating the problem as a mixed-integer second-order cone programming (MISOCP) problem, which can be solved via MPC with successive linearization [169]. In addition to thrusters, differential drag has also been used as the control for RPOs with possible fuel savings and without harmful jet firings [170]. By opening or closing the drag plates equipped with each spacecraft, the drag force acting on each vehicle can be modulated to control the RPO pro-cess. Different from the convex approach in [18], the differential-drag-based RPO has been formulated as a mixed-integer nonlinear programming (MINLP) problem because of the binary feasible control set and the free final time. By relaxing the control set from $\{-1, 0\}$ to $[-1, 0]$, a convex problem has been obtained, and its solution has proved to be a solution to the original problem, although a feasible control for the relaxed problem may not necessarily be feasible for the original problem [170]. The optimal final time was sought by solving a sequence of LP problems via a one-dimensional search.

In addition to fuel consumption, observability cost has also been considered and combined with fuel consumption to formulate a multi-objective, convex QP problem for bearings-only RPO missions [185]. Besides approaching a passive target, both the chaser and target can be controlled simultaneously for cooperative RPOs, which has been addressed via a combination of variable changes, constraint relaxation, and successive convexification [172]. The approach can be reduced into a single SOCP when linear dynamics are considered for linear impulsive RPO missions [186]. Additionally, mesh refinement and pseudospectral methods have been combined with convexification techniques for RPOs with improved solution accuracy while guaranteeing computational efficiency [173, 174]. Moreover, other relative dynamic models, such as those based on the Kustaanheimo-Stiefel transformation, have also been used to fundamentally shift the RPO problem formulation while still maintaining its adaptability to convex-optimization-based methods with potentially better solution accuracy [187]. Those dynamic models may be of particular interest to RPOs under non-impulsive low-thrust propulsion such as differential drag.

Furthermore, a number of publications have solved the RPO problems via QCQP and SDP perspectives. For example, using the popular linearized Tschauner-Hempel equations, impulsive control and continuous trajectory constraints have been merged to develop a convex description of the RPO problem with polynomial nonnegativity constraints, leading to an SDP problem (through a change of variables and introducing slack variables) that can be solved in polynomial time [175]. Inspired by the results of [175], an SDP-based glide-slope guidance algorithm has been proposed for minimum-fuel RPOs on elliptic orbits, and the formulation has been shown to be able to reduce into an LP problem when no trajectory constraints are enforced [176]. More recently, a QCQP-based method has been developed for multi-phase trajectory optimization problems with an application to a two-phase RPO via SDP relaxation and an SCP-type approach [188]. Different from the LMI-based or SDP-based approaches, an alternating minimization algorithm has been developed to solve a nonconvex QCQP transformed from a polynomial programming formulation of the RPO problem by solving a sequence of convex QP problems [189]. In the meantime, a multi-phase RPO problem has also been formulated as a general QCQP problem but solved using an alternating direction method of multipliers (ADMM) by introducing slack/auxiliary variables to simplify the sub-problems [190].

To address the iterative feasibility challenge as has been observed in many SCP-type algorithms, an iterative convex-optimization-based approach has been developed for soling nonconvex OCPs with linear dynamics and used to solve RPOs as a case study [191]. The approach is akin to SCP that requires solving a sequence of convex sub-problems; however, the method guarantees the feasibility of each intermediate iterate (given a feasible initial iterate) and facilitates monotonic convergence of the solution by formulating the original nonconvex problem as a difference of convex function (DC) programming problem [191]. More recently, a convex-concave decomposition method has been developed to address the shortcomings of the conventional linearization-based techniques for nonlinear equality constraints to facilitate the use of convex optimization for nonconvex OCPs [177]. By relaxing each nonlinear equality constraint into three convex or concave inequality constraints, the solution to the original problem can be obtained by solving a sequence of SOCP problems. This approach has shown to be effective in solving a 3-DoF fuel-optimal spacecraft circumnavigation problem where the deputy spacecraft is controlled into a specified relative orbit around the chief spacecraft [177]. In addition, a recent work in [178] has opened up another new perspective on how to solve nonconvex OCPs such as circumnavigation RPOs via convex optimization by showing that the original nonconvex problem can be cast into an equivalent two-variate constrained minimization problem that can be efficiently solved by a hybrid algorithm combining a convex relaxation method and a linearization-projection approach.

Other than 3-DoF cases, 6-DoF RPO problems have also addressed in the literature using convex optimization. For example, the method developed for 6-DoF PDG [105, 82] has proved to be effective in addressing 6-DoF RPOs by capturing the coupled translational and rotational dynamics using unit dual quaternions and formulating a convex QCQP problem to find the solution [192]. In addition, an SCP approach has been developed to solve a 6-DoF fuel-optimal RPO problem by transforming the original nonconvex problem into a series of SOCP sub-problems via successive convexification [179]. Specifically, the nonconvex field-of-view constraint was approximated as a second-order cone, while the concave obstacle-avoidance constraints were convexified into affine inequality constraints through linearization. In light of the effectiveness of state-triggered convexification techniques in solving PDG problems [106, 86, 110, 87], similar approaches have been developed to address 6-DoF fuel-optimal RPOs in the presence of mixed-integer constraints, such as plume impingement constraints that only need to be enforced when the two vehicles get close enough [180]. These constraints can be handled as state-triggered constraints within a continuous optimization framework via successive convexification without the need of solving difficult, time-consuming mixed-integer programming problems [181]. However, some unfavorable locking behavior has been observed in state-triggered SCP for RPO problems and can prevent the algorithm from converging [180]. More recently, a homotopy approach has been developed to address this phenomenon in solving 6-DoF RPOs with discrete logic constraints integrating numerical continuation and SCP into a single iterative solution process and approximating the discrete logic constraints with smooth functions using a homotopy parameter to control the approximation accuracy [193]. Of course, for simplicity, the 6-DoF RPO problem can also be decoupled into an attitude G&C problem and an orbit G&C problem, each of which can be solved by the existing convex optimization methods in the literature [194].

Similar to other applications, central to a successful G&C design for RPOs is the robustness of the method to uncertainties, disturbances, and maneuver execution errors, in particular when approaching a non-cooperative target, at a low computa-tional cost [195, 196]. Convex optimization has been combined with MPC to address these challenges by solving a convex optimization problem at each MPC iteration [197, 198, 199, 200]. Recently, a stochastic MPC method has been applied to solve RPOs in the presence of unbounded disturbances in [201], where the constraints such as obstacle avoidance were modeled as chance constraints that can be equivalently convexified into second-order cone constraints. In addition, covariance control has received a growing interesting in solving aerospace optimal G&C problems as discussed in subsection 3.1. Besides PDG, covariance control has also been used to solve RPO problems. For example, a nonlinear stochastic OCP has been formulated for a 6-DOF RPO trajectory optimization problem in a recent work considering initial state uncertainties and external disturbances as well as chance constraints on collision avoidance, sensor field-of-view, approach corridor, and control amplitude [182]. The stochastic problem was reformulated into a deterministic, convex form via relaxation, linearization, discretization, and introduction of auxiliary variables. An approximate optimal solution to the original problem was obtained by solving a series of SDP programs to simultaneously generate the nominal optimal trajectory along with affine feedback control policy [182].

### *3.4. Orbital Transfer*

A maneuver similar to RPOs is orbital transfer (OT), where the vehicle is guided along a transfer trajectory from its initial orbit to the destination in a target orbit for near-Earth orbital missions as well as interplanetary and deep space exploration missions (see Figure 5). Trajectory optimization for OTs has gathered increasing attention in the past two decades, spurred by novel propulsion technologies such as electric, nuclear, and solar-sail propulsion. Despite the highly appealing efficiency of novel propulsion systems, the produced thrust is usually very low, and the resulting trajectory optimization problem is challenging to solve. Different from short-range RPOs, long-range OTs may require many orbital revolutions when the initial and target orbits are widely spaced, and the long OT duration may lead to significant computational challenges [202]. Other major challenges associated with OTs include the high nonlinearity and severe coupling of state and control variables in the dynamics. In recent years, there has been a burgeoning development of highly efficient algorithms for potential real-time trajectory optimization and autonomous G&C of low-thrust OTs using convex optimization (see Table 4)

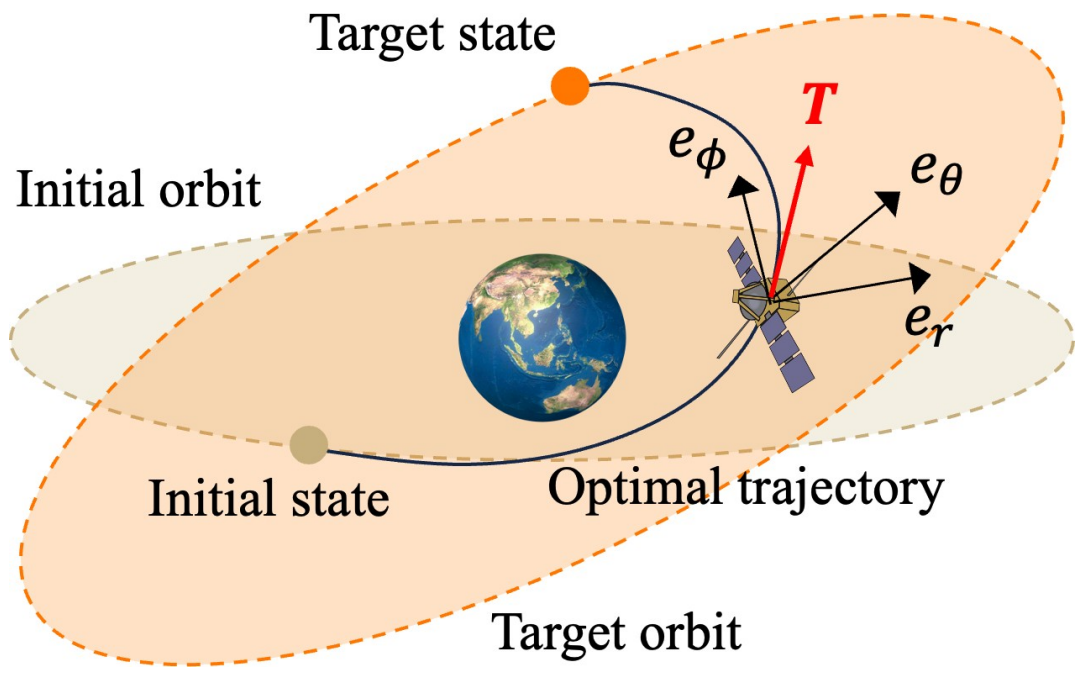

Figure 5: Schematic representation of orbital transfer (OT).

Table 4: Summary of representative publications on convex optimization for orbital transfer (OT).

| Reference | OT Problem | Approach | Slack Variables | Change of Variables | Convex Relaxation | Linear Approximation | Solver |
|---|---|---|---|---|---|---|---|
| [20] | 3-DoF fuel-optimal | SOCP + SCP | | ✓ | ✓ | ✓ | ECOS |
| [21] | 2-D time-optimal | SOCP + SCP | | ✓ | ✓ | ✓ | ECOS |
| [203] | 3-DoF time-optimal | SOCP + SCP | | ✓ | | ✓ | ECOS |
| [204, 205] | 3-DoF fuel-optimal in three-body systems | SOCP + SCP | | ✓ | ✓ | ✓ | MOSEK |
| [206] | 3-DoF fuel-optimal | SOCP + SCP + pseudospectral + indirect | | | ✓ | ✓ | ECOS |
| [207, 208] | 3-DoF fuel-optimal | SOCP + SCP + pseudospectral | | ✓ | ✓ | ✓ | ECOS |
| [209] | 3-DoF fuel-optimal multi-phase | SOCP + SCP + pseudospectral | | | ✓ | ✓ | ECOS |
| [210] | 3-DoF fuel-optimal | SOCP + SCP + homotopy | | ✓ | ✓ | ✓ | ECOS |
| [211, 212, 213, 214] | 3-DoF minimum closed-loop control | Covariance control + SOCP/SDP + SCP | | ✓ | ✓ | ✓ | MOSEK |

Earlier studies on convex-optimization-based OTs can be found in [215, 20, 21, 216], where both fuel-optimal and time-optimal low-thrust OTs have been solved with the aim of quickly obtaining optimal or near-optimal numerical solutions with high accuracy and low computational cost. To this end, the problems were formulated as general OCPs first, and then a series of transformation techniques was applied to convert the original problem into a convex formulation through a change of variables, relaxation of control constraints, and successive convexification and linear approximations. Specifically, inspired by the remarkable performance of SOCP for PDG and RPO problems [15, 16, 18], new state and control variables, such as $\tau = \frac{T}{m}, z = \ln m, \tau_r = \tau \cos \alpha_r$, $\tau_\theta = \tau \sin \alpha_r \sin \alpha_{\varphi\theta}$, and $\tau_\varphi = \tau \sin \alpha_r \cos \alpha_{\varphi\theta}$ were introduced to reduce the nonlinearity of the dynamic model. Lossless convexification was used to relax the nonconvex control constraint into a convex form, and an SCP algorithm was finally developed to find approximate optimal solutions to the original problem by solving a sequence of SOCP sub-problems. The equivalence of the relaxation and the existence of the solution to the relaxed problem have been proved [20, 21]. In addition, it is worth mentioning that a new independent variable with a monotonically increasing or decreasing trend may be needed to rewrite the equations of motion in a way that free-final-time problems can be readily transformed into convex optimization problems [21]. This is different from the popular approach where the free-final-time problem is converted into a fixed-final-time problem by normalizing the original time $t \in [t_0, t_f]$ to $\overline{t} \in [0, 1]$ with $t_f$ as an extra parameter to be optimized. This approach further increases the nonlinearity of the dynamic system because each dynamic equation must be multiplied by the $t_f$ parameter, which makes the dynamics more difficult to convexify.

With the initial success of convex optimization in solving OT problems in [20, 21], many similar approaches with various improvement mechanisms have emerged in addressing a variety of OT missions. For example, an SOCP-based SCP approach has been combined with the pseudospectral method to provide the adjoint variables via adjoint mapping to initialize a homotopic indirect method in solving a 3-DoF fuel-optimal low-thrust OT problem [206]. Comparison results have shown that the SCP method may suffer from a potential loss of solution optimality compared to the indirect method, although stable convergence of SCP is usually expected [217]. Also, the SOCP-based SCP method has been augmented by trust regions and virtual controls for improved convergence in solving a 3-DoF time-optimal solar-sail interplanetary trajectories by controlling the cone angle and clock angle of the sail [203]. Furthermore, SCP has been combined with sparse optimal control [218] to solve a 3-DoF fuel-optimal OT from a near rectilinear halo orbit to a low lunar orbit in the Earth-Moon system in the context of a circular restricted three-body problem (CR3BP) [219]. Later, the approach was employed to solve a low-thrust transfer between periodic halo orbits around the same libration point defined in a three-body problem [220, 204]. In the meantime, SCP has also been used to solve transfers between libration orbits in the Mars-Phobos system based on a CR3BP model [205].

In addition, convex optimization and convexification techniques have been applied to address a wide range of OT scenarios including ballistic capture [221], detection and estimation of spacecraft maneuvers [222, 223], OT optimization with variable specific impulse and engine shutdown constraint [224], multi-arc OT optimization with constraints on duration of arcs and linkage constraints [209], multi-phase gravity assist low-thrust trajectory optimization under multi-body dynamics [225], and space intercept with nonlinear terminal constraints [226]. Techniques, such as homotopic approaches [227, 210], mesh refinement [228], Radau pseudospectral [229, 230, 207, 231], Chebyshev pseudospectral [232], and differential-algebra-based trust regions [233], have been employed to enhance the performance of SCP in solving OT problems. More recently, the performance of different trust-region methods (hard/soft trust region with constant/varying thrust-region shrinking) and discretization approaches (adaptive Legendre–Gauss–Radau pseudospectral methods, an arbitrary-order Legendre– Gauss–Lobatto technique based on Hermite interpolation, and a first-order-hold interpolation method) has been evaluated on solving low-thrust OT problems using SCP [208]. In addition, the impacts of different coordinate representations have been assessed in SCP-based OTs. The modified orbital/equinoctial elements, spherical, and cylindrical coordinates seem to outperform Cartesian coordinates in terms of success rate [234].

When it comes to handling uncertainties and disturbances, a convex-concave procedure has been used to convert the original, nonconvex OT problem with chance constraints into a sequence of convex sub-problems for risk-aware trajectory design by quantifying the uncertainties of orbital states [235]. Recently, chance-constrained covariance control has been used to formulate the OT problem as a stochastic OCP where the vehicle dynamics is modeled as a stochastic system that is steered from an initial probability distribution to a desired probability distribution subject to probabilistic state and control constraints modeled as chance constraints [211, 213]. Through a change of variables and constraint relaxation, the covariance matrix propagation was transformed into a set of semidefinite cone constraints and the original covariance control problem was reformulated as an SDP problem with proved lossless relaxation property [236, 212]. Finally, an SDP-based SCP approach was established and used to simultaneously generate the optimal nominal transfer trajectory and the feedback control policy to compensate the flight uncertainties and disturbances. More recently, the approach has been extended to address low-thrust OTs accounting for mass uncertainty where the propagation of the mean and covariance of mass is approximated by a set of convex constraints via a change of variables [214].

### *3.5. Spacecraft Constrained Reorientation*

Many space missions require the spacecraft to change its orientation for specific mission purposes in the presence of attitude constraints. For example, the direction of the spacecraft's high-gain antenna may need to remain in a particular cone for communication with ground stations during the reorientation. In addition, some sensitive onboard instruments (e.g., infrared telescopes, interferometers, star trackers) may be kept away from direct exposure to bright objects such as the sun during the attitude maneuver. A key component for such operations is the G&C algorithm that can be run in an autonomous manner onboard to produce appropriate steering laws for safe and efficient reorientation maneuvers. The problem can be formulated as an OCP that finds the control torques that optimize an objective function over a time interval subject to the initial and final states, nonlinear attitude kinematics and dynamics, bounded angular velocities and control inputs, and constraints on the orientation of the spacecraft (see Figure 6). This problem is referred to as space-craft constrained reorientation (SCR) in this paper. The difficulties in solving this problem are attributed to the nonlinear attitude dynamics and nonconvex attitude constraints. The following reviews the representative works on optimal SCR via convex optimization (see Table 5).

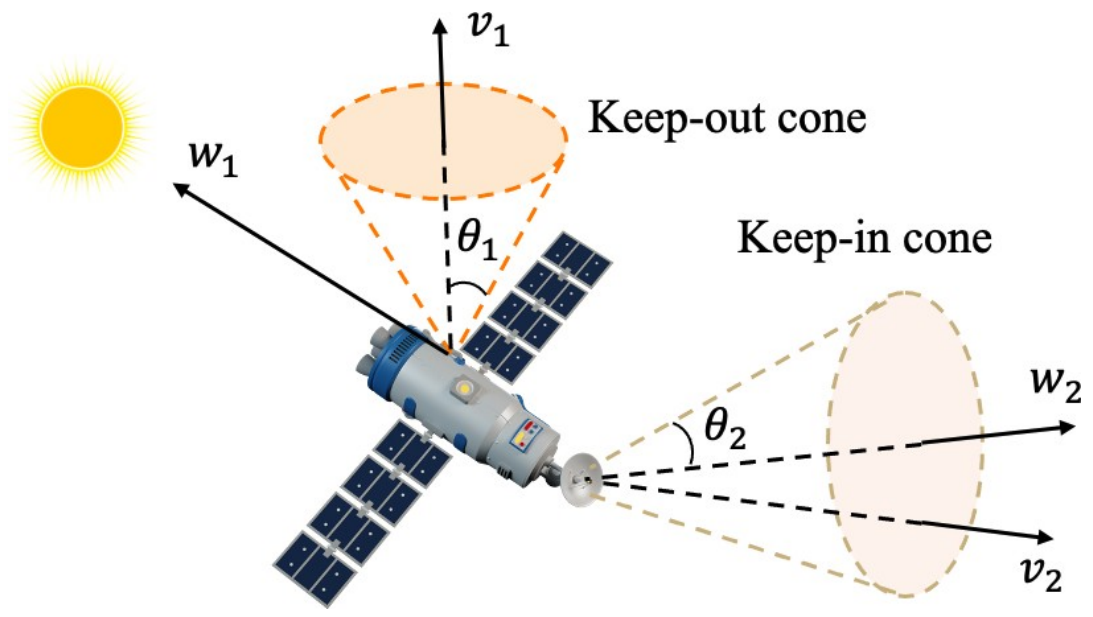

Figure 6: Schematic representation of spacecraft constrained reorientation (SCR).

Table 5: Summary of representative publications on convex optimization for spacecraft constrained reorientation (SCR).

| Reference | SCR Problem | Approach | Slack Variables | Change of Variables | Convex Relaxation | Linear Approximation | Solver |
|---|---|---|---|---|---|---|---|
| [237, 238, 239, 240] | Minimum terminal error | SDP + SCP | | | ✓ | ✓ | - |
| [241] | Multi-objective | SDP + MPC | | | ✓ | ✓ | SeDuMi |
| [242] | Minimum angular velocity and angular acceleration | MICP | ✓ | | ✓ | | Gurobi |
| [243] | Multi-objective | MICP | ✓ | | ✓ | ✓ | Gurobi |
| [244] | Minimum control | QCQP + SDP + SCP | ✓ | | ✓ | | - |
| [245] | Multi-objective | SDP + SCP | ✓ | | ✓ | ✓ | MOSEK |
| [246] | Minimum energy | QP + SCP + line search | | | | ✓ | ECOS |
| [247] | Minimum energy | SOCP + SCP + pseudospectral | | | ✓ | ✓ | MOSEK |
| [248] | Minimum power consumption | SOCP + SCP | ✓ | | | ✓ | ECOS |

Motivated by the advances in LMI theory for solving optimization problems defined over matrix spaces, the SCR problem has been approached via an SDP-based strategy by exploiting the nonconvexity of attitude constraints as well as the nonlinearity of system dynamics [237, 238]. The challenge lies in how to deal with the nonconvex quadratic constraints on the orientation of the spacecraft in the form of $v^T w \leq \cos\theta$, where $v$ and $w$ represent the unit vectors describing the boresight direction of sensitive onboard instruments and the direction of the undesired celestial object to be avoided in an inertial coordinate frame, respectively. The angle $\theta \in [0, \pi]$ defines the required minimum angular separation of these two vectors. To facilitate the implementation of efficient convex optimization methods, a quaternion representation of the spacecraft attitude has been used, leading to an equivalent quaternion characterization of the attitude constraints like $q(t)^T \widetilde{A} q(t) \leq 0$, which is nonconvex in the spacecraft orientation since $\widetilde{A}$ is not positive semidefinite. A key step is to relax this nonconvex constraint into an equivalent convex quadratic inequality or an LMI [237]. Combined with a linear approximation to the dynamic equations, a solution to the original SCR problem can be obtained by iteratively solving an SDP problem. This approach has been augmented to addressed different types of attitude constraints (hard or soft, static or dynamic) [238].

The underlying ideas in [237, 238] has contributed to the development of many convex optimization methods for SCR problems. For example, the SDP relaxation approach has been implemented in an MPC framework via linearization of the spacecraft attitude dynamics under similar attitude constraints [241]. Also, the results have been expanded to develop a potential-function-based approach to SCR with different types of attitude-constrained (forbidden or mandatory) zones defined by unit quaternions [239, 240]. A convex parameterization of these zones has been utilized to construct a logarithmic barrier potential for the synthesis of smooth and strictly convex attitude control laws. Interestingly, a recent work has pointed out that this approach may suffer from “undesired equilibria” by showing that the relaxed keep-out zones are actually non convex and the designed potential barriers are not convex functions, which needs further investigation [249]. Regardless of this controversy, an SCR problem subject to pointing and angular rate constraints has been solved by mixed-integer convex programming (MICP) by leveraging lossless convexification of nonconvex quadratic pointing constraints in [238, 240] and using binary variables to enforce the unity constraint on the quaternion magnitude [242]. MICP has also been used to solve SCR with both keep-in (inclusion) and keep-out (exclusion) pointing constraints by introducing binary variables to the formulation in [237, 238] to define logical pointing constraints due to redundant sensors and relaxing the set of nonconvex quadratic attitude constraints into mixed-integer convex constraints [243]. Meanwhile, SCR under similar constrained zones has been formulated as a general nonconvex QCQP and relaxed into an SDP with rank one constraint [244]. An iterative rank minimization approach was development to find this rank one matrix and converge to an optimal solution. In addition, [245] expanded on the results in [238, 240] and formulated the SCR problem as an SDP using the direction cosine matrix directly with reaction wheels for controls.

Rather than using the SDP relaxation and solving a sequence of SDP problems, a minimum-energy fixed-time rest-to-rest SCR problem has been addressed via a QP-based SCP method for an asymmetric rigid-body spacecraft in [246]. Through successive convexification, the solution was sought by solving a sequence of convex QP problems, which aids computational efficiency due to QP's lower complexity and availability of more applicable solvers than SDP. A line search was also introduced to promote the convergence of the SCP method that has been shown to converge even with trivial initial trajectories [246]. With the same goal of improving the computational efficiency of the algorithm, a set of convexification techniques has been combined with the Gauss pseudospectral method to relax an energy-optimal SCR problem into a series of SOCP problems [247]. More recently, an SOCP-based SCP has been used to solve a minimum-power-consumption SCR problem onboard the Satellite for Optimal Control and Imaging (SOC-i) CubeSat as part of its G&C flight software [248].

### *3.6. Space Robotic Manipulation*

The emerging active debris removal (ADR) and on-orbit servicing, assembly, and manufacturing (OSAM) technologies require a space robotic system (i.e., a base space-craft equipped with one or more manipulators) to support a variety of robotic missions such as capturing tumbling space objects, building large space structures, and refueling or fixing on-orbit satellites (see Figure 7) [250]. However, the mobile base spacecraft platform, highly nonlinear coupled base-manipulator dynamics, and complex operation constraints pose significant G&C challenges for such maneuvers. Also, the possibly unknown properties (e.g., mass, moment of inertia, shape) and motion characteristics (e.g., rotational rate) of the target (e.g., a tumbling space debris or a malfunctioning satellite) add another level of mission complexity and require performing target identification, motion prediction, and real-time decision-making autonomously onboard. Inspired by the remarkable success of convex-optimization-based G&C techniques in the areas of PDG and RPO as reviewed above, multiple approaches have been proposed to convexify the OCP formulation of the space robotic manipulation (SRM) problem for potential real-time onboard applications (see Table 6).

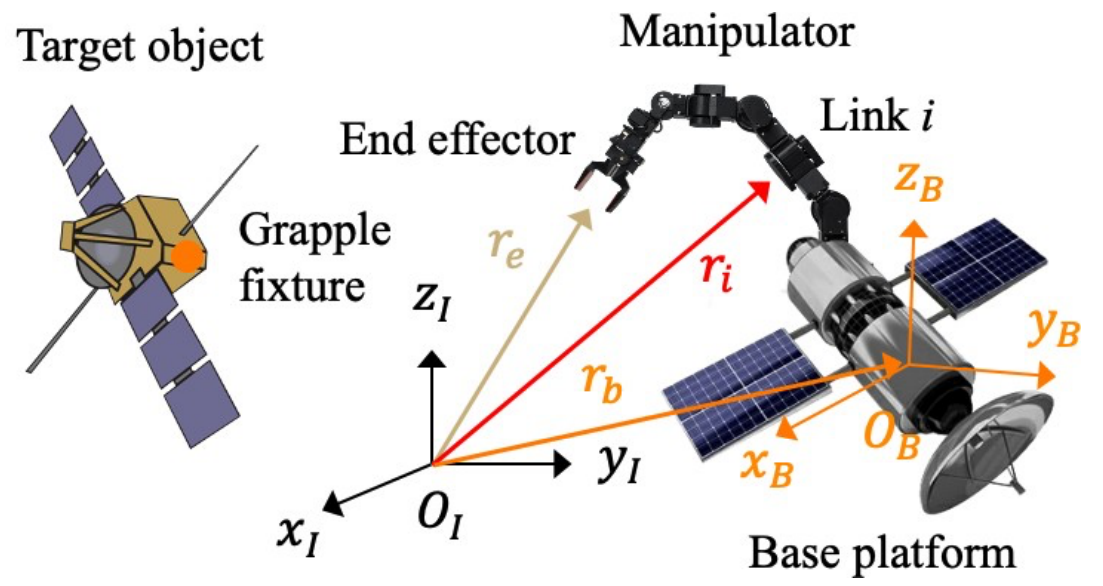

Figure 7: Schematic representation of space robotic manipulation (SRM).

Table 6: Summary of representative publications on convex optimization for space robotic manipulation (SRM).

| Reference | SRM Problem | Approach | Slack Variables | Change of Variables | Convex Relaxation | Linear Approximation | Solver |
|---|---|---|---|---|---|---|---|
| [251, 252, 253] | 3-DoF minimum control | QP + SCP | | | | ✓ | SDPT3 |
| [22] | 10-DoF minimum control and base attitude | QP | ✓ | | | ✓ | Gurobi |
| [254] | 8-DoF minimum control and kinetic energy | QP + SCP | | | | ✓ | Gurobi |
| [255] | 10-DoF minimum control and manipulability | QP | | | | ✓ | CSDP |
| [256] | 7-DoF minimum contact forces | SDP + SCP | | | ✓ | ✓ | - |
| [257] | 7-DoF minimum control and base attitude | QCQP + SCP + pseudospectral | | | | ✓ | MOSEK |
| [258] | 6-DoF minimum time | SOCP | ✓ | | | ✓ | SeDuMi |

One of the earliest works on using convex optimization for SRM G&C was reported in [251], where the SRM maneuver was divided into two sub-tasks. The first task aims to solve the system-wide center-of-mass translation problem that generates the control profile required to translate the base spacecraft to a location close enough to the target object. The second task solves the internal re-configuration problem that produces the control histories to re-orient the base spacecraft and re-configure the manipulator. Both tasks were carried out simultaneously and solved individually by SCP approaches. The inherent nonlinear dynamics and nonconvex constraints (e.g., collision avoidance and line-of-sight) were handled by successive convexification [251]. Dividing the SRM process into sub-maneuvers simplifies the optimization and operation to some extent at the cost of degrading the optimality of the overall solution. More theoretical convergence analysis and hardware-in-the-loop experiments have been presented to validate these SCP-based SRM techniques [252, 253].

In the meantime, a convex QP approach has been developed for trajectory planning of redundant manipulators on a free-floating mobile spacecraft platform (no base actuation) with nonzero initial momentum subject to bounded joint angles, joint velocities, and joint accelerations as well as obstacle avoidance and end-effector constraints [22]. Great effort was devoted to obtaining a convex QP formulation with linear constraints through relaxation of the nonlinear equality end-effector pose constraints and the nonconvex obstacle avoidance constraints. By solving a resulting convex QP on each discrete node, optimal collision-free end-effector trajectories can be generated to minimize the base attitude disturbance and control effort while satisfying the joint limits and end-effector task constraints [22]. Shortly after, the QP-based approach was combined with SCP to solve a point-to-point SRM planning problem, aiming to find an optimal joint path, along which the manipulator drives the end effector from its initial state to a desired target condition under obstacle avoidance and end-effector pose constraints [254].

In addition, a convex QP approach has been developed to solve a trajectory planning problem for a free-floating space manipulator incorporating constraints on end-effector trajectory tracking and spacecraft attitude stabilization as well as the joint angle/velocity/acceleration constraints [255]. More recently, an SDP approach has been introduced to solve a post-capture grasping force optimization problem for a dual-arm spacecraft [256]. The original nonconvex problem was relaxed into an SDP problem subject to base force/torque limits, joint torque limits, and LMIs resulted from converting the nonlinear friction constraints into the positive definiteness of specific symmetric or Hermitian matrices. Minimum-contact-forces solutions were obtained by iteratively solving the relaxed SDP problem with the aid of a line search method. Further, a Legendre pseudospectral method has been used to improve the efficiency of the SCP approach in solving a point-to-point trajectory planning problem for a free-floating space robotic system with a 7-DoF manipulator [257]. Most recently, a time-optimal path tracking problem with dynamic and base velocity constraints has been addressed for a 6-DoF dual-arm free-floating space manipulator and transcribed into an SOCP problem through introduction of slack variables and affine approximations [258].

### *3.7. Spacecraft Formation Flying and Station Keeping*

Coordinating a group of smaller distributed spacecraft in formations or configurations (see Figure 8) can accomplish some space missions that are difficult or impossible for a larger, more expensive, monolithic spacecraft. Spacecraft formation flying (SFF) can bring significant benefits over single vehicles including higher redundancy, simpler designs, faster response times, and cheaper replacement. These traits make SFF ideal for a variety of missions such as reconnaissance, observation, communication, meteorology, and terrain mapping [259]. Orbital station-keeping (SK), a series of active orbital maneuvers that compensate for orbital perturbations, is vital for both single-spacecraft and multi-spacecraft scenarios to maintain a stationary orbit or configuration. This is of particular importance for spacecraft in a halo orbit around a libration point where the orbit is unstable and obvious deviations in position and velocity may occur if no active control is employed [260]. To achieve the goals of SFF and SK, a flexible, robust, and computationally efficient G&C framework that can be used as part of the onboard system is needed to plan optimal maneuvers while satisfying constraints such as collision avoidance and plume impingement. Several convex-optimization-based approaches have been explored in the literature for SFF and SK problems. Some of the representative research efforts are summarized in Table 7 and briefly reviewed below.

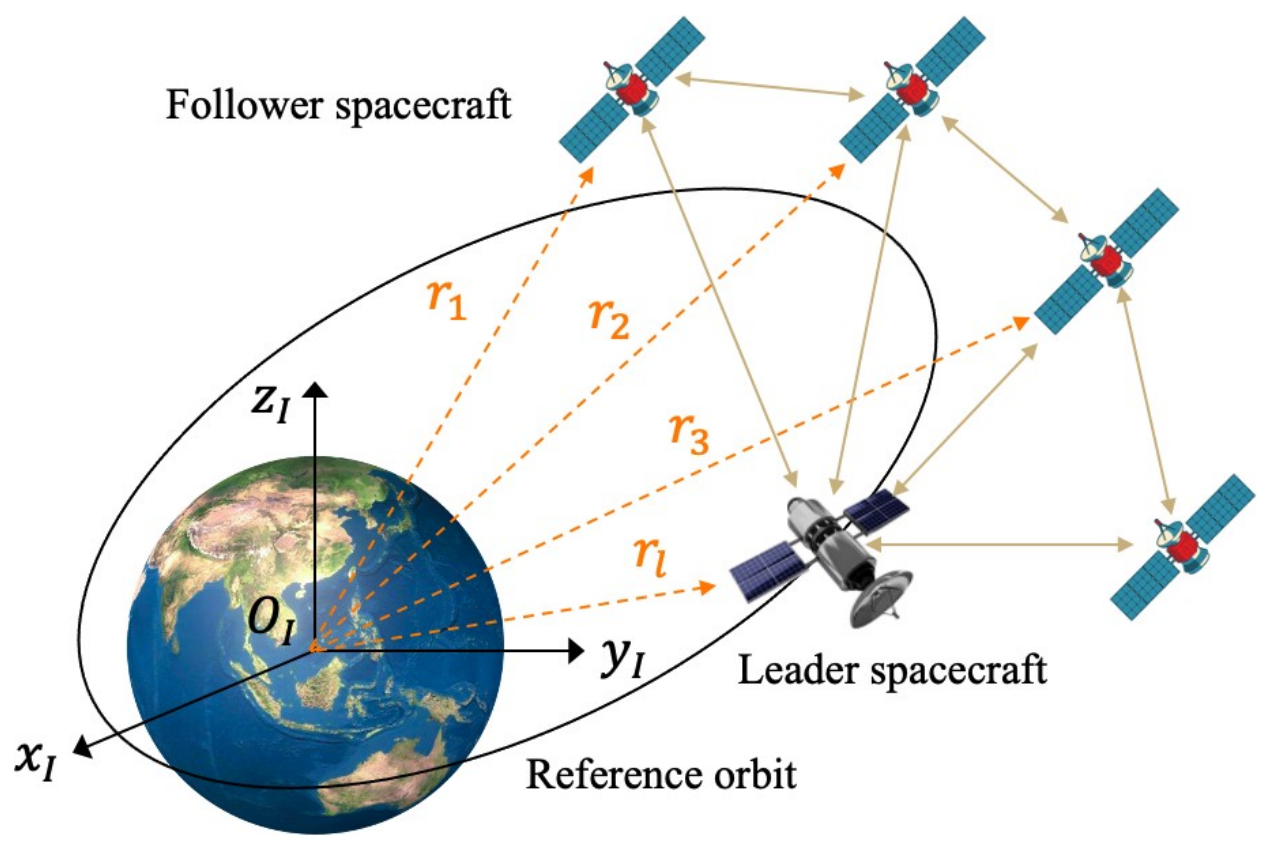

Figure 8: Schematic representation of spacecraft formation flying (SFF).

Table 7: Summary of representative publications on convex optimization for spacecraft formation flying (SFF) and station keeping (SK).

| Reference | Problem | Approach | Slack Variables | Change of Variables | Convex Relaxation | Linear Approximation | Solver |
|---|---|---|---|---|---|---|---|
| [149, 261, 150] | 3-DoF fuel-optimal SFF | LP/MILP | ✓ | | | ✓ | CPLEX |
| [262] | 3-DoF fuel-optimal SFF | SDP/SOCP | ✓ | ✓ | ✓ | ✓ | SeDuMi |
| [263] | 3-DoF minimum control SFF | SDP + SCP | | | ✓ | ✓ | - |
| [264, 265] | 3-DoF fuel-optimal SFF | SOCP + SCP | | | ✓ | ✓ | SDPT3, MOSEK |
| [266] | 3-DoF multi-objective SFF | SOCP + SCP | ✓ | | ✓ | ✓ | ECOS, MOSEK |
| [267] | 3-DoF fuel-optimal SK | LP | ✓ | | ✓ | ✓ | - |
| [268, 269] | 3-DoF fuel-optimal SK | LP | ✓ | | | ✓ | MOSEK |
| [270] | 3-DoF minimum control and tracking error SK | SDP + MPC | ✓ | | ✓ | | MOSEK |
| [271] | 3-DoF minimum control and tracking error SK | SOCP + MPC | ✓ | | ✓ | ✓ | CVX |

Earlier works have shown the efficiency of LP and MILP in addressing SFF problems with linearized orbital dynamics and mixed linear/integer constraints for collision avoidance and plume avoidance [149, 261, 150]. Other than simple linear approximations, relaxation techniques such as introduction of slack variables and change of variables have also been used to transform the nonconvex SFF G&C problem into more efficient convex forms such as SDP or SOCP [262]. In the presence of attitude forbidden zones, the quadratic convex constrained zone formulation in [239, 240] has been utilized for a group of spacecraft to achieve identical orientation [272]. In addition to solving single convex problems, an iterative SDP relaxation method has been used to handle SFF problems, aiming to determine a formation and topology of a group of spacecraft modules with guarantees on network connectivity while minimizing the total control effort [263]. SCP approaches have also been developed for SFF to produce collision-free fuel-efficient reconfiguration trajectories of spacecraft swarms through successive linearization and convex relaxation of the nonconvex collision-avoidance constraints [273, 264, 274, 265]. These approaches can be implemented in a centralized or a decentralized manner [275, 276] for small-scale or large-scale formations [277, 278] with free-flying or tethered configurations [279]. The SCP approach has been combined with Markov Chains, MPC, and pseudospectral method to address collision-free formation flying of large-scale spacecraft swarms [280, 281]. More recently, the SCP method has been improved with appealing convergence to solve optimal reacquisition planning problems for distributed spacecraft systems in the context of gravitational wave detection [266].

SK constraints have also been considered in solving optimal G&C problems in the context of multi-agent space missions such as SFF. For example, an avoidance planning problem considering SK constraints has been transformed into and solved as an LP problem (no binary variables) for a distributed set of close spacecraft to produce fuel-efficient maneuvers while maintaining the desired station on orbit [267]. Inspired by [261], station-keeping maneuvers have been determined by formulating and solving an LP problem for a geostationary satellite through a novel affine formulation of the equations of motion in [268], which has been improved to determine station-keeping maneuvers for a fleet of satellites in a geostationary slot by explicitly considering the thruster configuration and incorporating each individual thruster control in the problem formulation [269]. Via a different approach, [270] posed the SK problem as an SDP problem based on a polynomial approximation of the nonlinear CR3BP model to generate the optimal strategies for SK of halo orbits at $L_1$ libration point for the Sun-Earth three body system. More recently, an SOCP-based approach has been developed and implemented under MPC framework to address the SK control problem of halo orbit in the Earth-Moon CR3BP system via linearization of the dynamics and convexification of the nonconvex control constraints [271]. In addition, the sparse optimal control technique from [219, 218] has been used to solve an SK problem around libration point orbits in a Sun-Earth CR3BP system by formulating a convex optimization problem based on the Floquet-Lyapunov transformation of the dynamics [282].

## 4. Applications to Air Vehicles

In addition to space systems, convex optimization has also found many applications in the development of optimal G&C methods for atmospheric flight vehicles including hypersonic/entry vehicles, missiles and projectiles, launch/ascent vehicles, and low-speed manned/unmanned air vehicles, which will be surveyed in this section.

### *4.1. Hypersonic/Entry Guidance*

Hypersonic flight has been a critical phase for many space missions such as Earth reentry (e.g., Space Shuttle), planetary entry (e.g., Mars entry), and hypersonic weapons. Due to the high-speed atmospheric flight, the main purpose of hypersonic G&C is to control the variation (usually dissipation) of the vehicle's kinetic energy to meet specific mission requirements while satisfying various constraints (see Figure 9). Closely related to planetary entry/reentry, aero-assisted maneuvers also experience hypersonic atmospheric flight to either capture the vehicle into a closed orbit around the target planet (i.e., aerocapture) or achieve a large change in the direction of the velocity (i.e., aerogravity assist) by controlling the aerodynamic forces for reduced propulsion requirements. However, hypersonic/entry guidance (HEG) problems are difficult to solve due to the highly nonlinear dynamics, nonconvex path constraints (e.g., heat rate, normal load, dynamic pressure), and possible waypoint and no-fly zone constraints. Convex-optimization-based methods have received significant attention (see Table 8) for potential real-time hypersonic/entry trajectory generation and autonomous G&C due to their fast computational speed, easy implementation, and ability to enforce common constraints for various types of HEG missions.

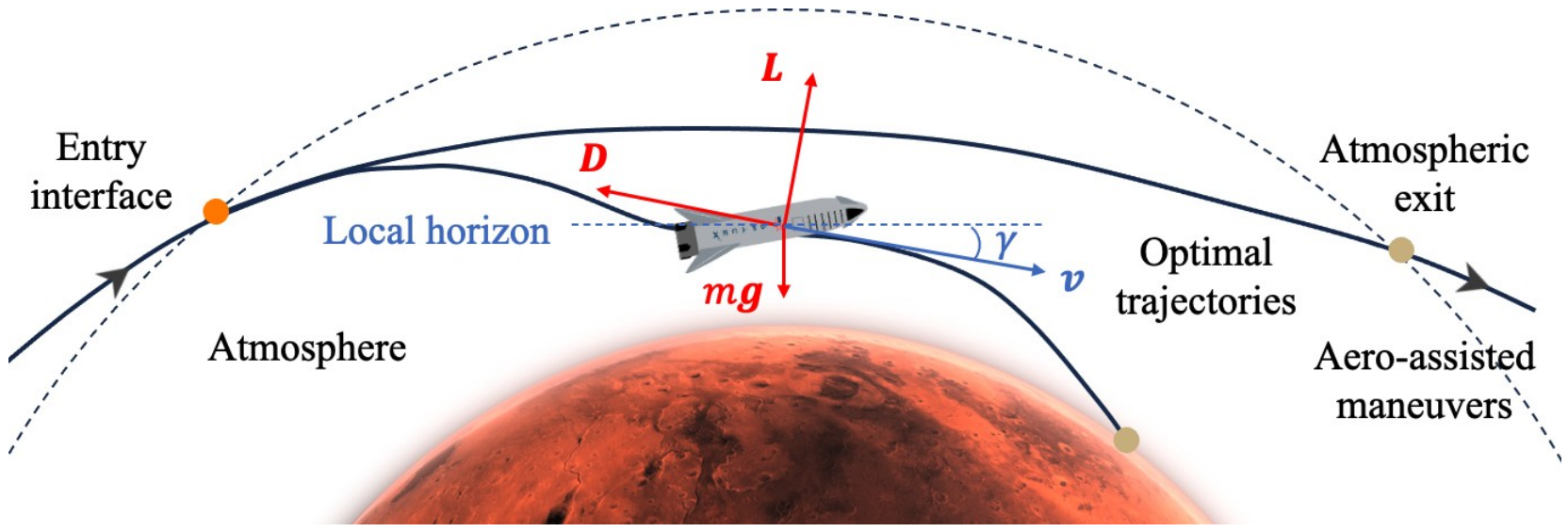

Figure 9: Schematic representation of hypersonic/entry guidance (HEG).

Table 8: Summary of representative publications on convex optimization for hypersonic/entry guidance (HEG).

| Reference | HEG Problem | Approach | Slack Variables | Change of Variables | Convex Relaxation | Linear Approximation | Solver |
|---|---|---|---|---|---|---|---|
| [25, 283, 24] | 3-DoF minimum-time, minimum-heat-load, minimum-oscillation, and maximum-crossrange | SOCP + SCP | ✓ | ✓ | ✓ | ✓ | MOSEK |
| [26] | 3-DoF minimum-terminal-velocity and minimum-heat-load | SOCP + SCP | | ✓ | | ✓ | SDPT3 |
| [27] | 3-DoF maximum-terminal-velocity and minimum-heat-load | SOCP + SCP + line search and trust region | ✓ | ✓ | | ✓ | ECOS |
| [88] | 3-DoF minimum- and maximum-peak-normal-load | MICP + SCP + line search | ✓ | ✓ | | ✓ | ECOS, Gurobi |
| [284] | 3-DoF minimum-time multi-phase | SOCP + SCP | ✓ | ✓ | ✓ | ✓ | ECOS |
| [285] | 3-DoF minimum-time | SOCP + SCP + pseudospectral | | ✓ | ✓ | ✓ | ECOS |
| [286] | 3-DoF minimum-time | SOCP + SCP + pseudospectral + trust region | | ✓ | | ✓ | MOSEK |
| [287] | 3-DoF maximum-impact-velocity | SOCP + SCP + line search | ✓ | ✓ | ✓ | ✓ | MOSEK |
| [288] | 3-DoF minimum-time and maximum-crossrange | Chance-constrained SCP | | | ✓ | ✓ | - |
| [28] | 3-DoF minimum-curvature | SOCP + pseudospectral | ✓ | ✓ | ✓ | ✓ | ECOS |
| [289, 290] | 3-DoF closed-loop tracking guidance | QCQP + pseudospectral | | ✓ | | ✓ | ECOS, MOSEK |
| [291] | 3-DoF minimum-impulse, minimum-time, and minimum-heat-load aerocapture | SOCP + SCP | | ✓ | ✓ | ✓ | MOSEK |
| [292] | 3-DoF minimum-$\Delta V$ aerocapture | SCP + covariance control | | ✓ | | ✓ | - |
| [293] | 2-D minimum-$\Delta V$ aerogravity assist | SOCP + SCP | | ✓ | ✓ | ✓ | MOSEK |

Earlier publications on convex-optimization-based HEG focused on the challenges of relaxing nonlinear and nonconvex control terms in the flight dynamics into convex forms that can be handled by IPMs [25, 283, 24]. Assuming a predetermined velocity-dependent angle-of-attack profile based on thermal protection and range considerations, the bank-angle components in the entry dynamics can be replaced by new controls, e.g., $u_1 = \cos\sigma$ and $u_2 = \sin\sigma$. Great efforts have been devoted to demonstrating the equivalence of relaxing the equality constraint $u_1^2 + u_2^2 = 1$ into the inequality constraint $u_1^2 + u_2^2 \leq 1$, i.e., assuring that the optimal solution of the relaxed problem lies on the boundary of the control set. Combining this relaxation technique with successive linearization and using the energy-based equations of motion, a successive SOCP approach has been developed to solve minimum-time [25], minimum-heat-load [25], minimum-oscillation [283], and maximum-crossrange [24] entry problems.

In the meantime, the time-based equations of motion have also been used to formulate and solve HEG problems. To avoid nonconvex control constraints and facilitate potentially more accurate solutions, the equations of motion can be reformulated by defining bank-angle rate as the new control with an additional state equation $\dot{\sigma} = u$. Through successive convexification, the minimum-terminal-velocity and minimum-heat-load problems have been solved via an SOCP-based SCP approach [26]. Aiming to improve the convergence of the algorithm, the approach was later improved by the line-search and trust-region techniques for HEG problems [27], including the minimum- and maximum-peak-normal-load hypersonic/entry trajectory optimization [294, 88]. The maximum-peak-normal-load problem was posed as a discrete-event max-max OCP, which was transformed into a sequence of MICP problems through a combination of a Big-M method and a line-search SCP [88].

In the past years, the SCP approach has been improved with the aid of techniques such as pseudospectral methods [285, 286], adaptive mesh refinement [295], *hp*-adaptive pseudospectral discretization [296], and virtual control [297] and have been applied for solving a wide range of HEG problems such as maximum-impact-velocity spiral-diving trajectories [287], multi-phase missions [284, 298], high-accuracy

HEG trajectory optimization with no-fly zones [299], and trajectory optimization under probabilistic constraints [288, 292]. Specifically, in the presence of waypoint constraints, the problem has been solved via a multi-phase SOCP-based SCP approach [284]. The reentry and landing phases have also been combined and solved as a multi-phase problem, which has been addressed via pseudospectral SCP [298]. In the presence of uncertain constraints, the hypersonic trajectory optimization problem has been formulated as a stochastic OCP where the probabilistic constraints are modeled as chance constraints [300]. Through a smooth and differentiable approximation of the probabilistic constraints, the original chance-constrained stochastic OCP can be transformed into a deterministic version that has been solved by an SCP approach [288]. Among the improved SCP approaches, the trust-region-based SCP has been investigated in a recent publication [301], showing that the trust-region order has obvious effects on the optimality of the converged solution, and higher-order trust-region SCP algorithms have been shown to outperform lower-order ones using HEG trajectory optimization as a case study [301].

Other than the successive approaches, a single SOCP problem has been formulated based on drag-energy dynamics and solved in combination with pseudospectral method to generate feasible drag-energy profiles through introduction of a set of new variables for the inverse of the drag acceleration [28]. Also, in addition to SOCP-based approaches, a sequential SDP approach has been explored for HEG problems by formulating the problem as a polynomial OCP and then a general QCQP via introducing new variables and quadratic constraints. An SDP relaxation technique has been utilized to relax the nonconvex QCQP problem into a sequence of SDP problems [302]. The convergence of this successive approach may be proved; however, the study of this approach has not been continued because solving SDP problems is generally much more time-consuming than solving SOCP problems. Further, besides generating optimal reference hypersonic/entry trajectories, convex optimization has been used for closed-loop tracking guidance for hypersonic/entry vehicles. For example, a convex QCQP problem has been solved in each guidance cycle to track the optimal reference trajectories generated by the successive SOCP approach [303, 289]. This numerical closed-loop HEG approach has recently been enhanced by using the Legendre–Gauss–Radau pseudospectral method for discretization to improve the computational efficiency while preserving solution accuracy for Mars entry [290]. Also recently, convex optimization has been combined with the popular predictor-corrector guidance algorithm for robust HEG by solving a single convex trajectory optimization problem for correction plan [304, 305].

In addition to entry and reentry, convex optimization has been used to develop real-time G&C algorithms for the atmospheric flight portion of aerocapture and aero-gravity assist (AGA) maneuvers. For example, successive convexification has been used to minimize the $\Delta V$ correction for the ADEPT (Adaptable Deployable Entry Placement Technology) planetary entry vehicle through active bank angle modulation for Mars aerocapture missions considering nonlinear dynamics and nonlinear boundary conditions [23]. The convex relaxation techniques in [25] have been extended to develop an SOCP-based SCP algorithm for a series of optimal aerocapture problems including minimum-impulse, minimum-time, and minimum-heat-load problems [291]. To explicitly consider model uncertainties for aerocapture G&C, chance-constrained covariance steering has been applied to jointly optimize updates to the feedforward control inputs and the corresponding feedback gains via an SCP approach [292]. Besides capturing the spacecraft into an orbit around the target planet, AGA maneuvers, using aerodynamic forces to augment gravity and achieve a larger change in direction than aerocapture, have also been solved through SCP [293].

### *4.2. Missile/Projectile Guidance*

Guiding an aerodynamically controlled missile or projectile to impact a stationary or mobile target has received sustained attention for decades. To achieve the best warhead effectiveness, the missile/projectile is expected to hit the target as accurately as possible along a specific direction (i.e., impact angle) as illustrated in Figure 10. In addition to the miss distance and impact-angle constraints, advanced missile/projectile guidance (MPG) should also consider control limits, field-of-view constraints, and possible constraints on dynamic pressure and heat rate [306, 307]. It is vital yet challenging to generate feasible and even optimal MPG commands and corresponding trajectories for different mission scenarios considering state, control, impact angle, and various path constraints. Convex optimization provides an efficient numerical approach to addressing such complicated problems with low computational cost and high solution reliability. In the atmospheric flight vehicle domain, we found a few publications on using convex optimization for MPG problems (see Table 9), which are briefly reviewed in this subsection below.

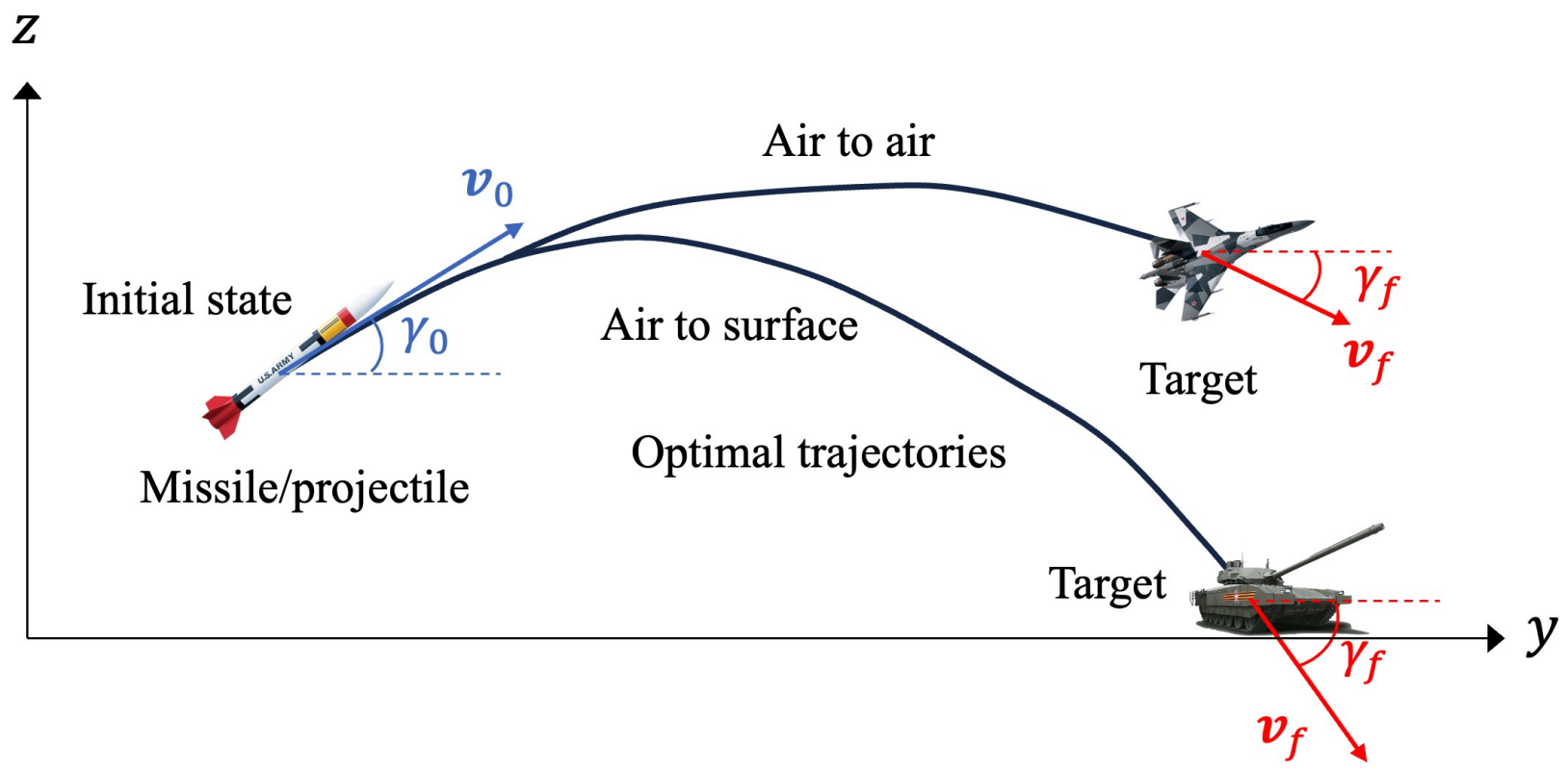

Figure 10: Schematic representation of planar missile/projectile guidance (MPG).

Table 9: Summary of representative publications on convex optimization for missile/projectile guidance (MPG).

| Reference | MPG Problem | Approach | Slack Variables | Change of Variables | Convex Relaxation | Linear Approximation | Solver |
|---|---|---|---|---|---|---|---|
| [308] | 3-DoF maximum-impact-velocity | SOCP + SCP | | ✓ | ✓ | ✓ | MOSEK |
| [309] | 2-D minimum-control | SOCP + SCP | ✓ | ✓ | | ✓ | MOSEK |
| [310] | 2-D maximum-impact-velocity | SOCP + SCP | | ✓ | ✓ | ✓ | MOSEK |
| [311] | 2-D minimum-energy | SOCP + SCP | ✓ | ✓ | | ✓ | MOSEK |
| [312] | 2-D minimum-energy | SOCP + SCP + L1 penalty | ✓ | ✓ | ✓ | ✓ | MOSEK |
| [313] | 3-DoF minimum-time | SOCP + SCP + pseudospectral | | | ✓ | ✓ | CVX |
| [314] | 3-DoF minimum-fuel and minimum-time | QP + SCP + pseudospectral | ✓ | | | ✓ | - |
| [315] | 2-D maximum-impact-velocity multi-phase | SOCP + SCP + pseudospectral + thrust region | ✓ | | | ✓ | MOSEK |
| [316] | 3-DoF minimum-control | SOCP + MPC | | | | ✓ | MOSEK |

The earliest attempt observed in this area was a trajectory optimization problem solved for an aerodynamically controlled missile to impact a ground target via an exact convex relaxation approach [308], where the terminal flight phase of the missile was optimized by solving a sequence of SOCP problems with both angle of attack and bank angle as the controls under impact angle and dynamic pressure constraints. An immediate issue faced by this approach was the exactness of the relaxation technique after introducing a set of new control variables. As observed in other applications, the relaxation may not be exact in general when the state inequality constraint becomes active. This issue was hurdled by introducing a regularization term, and theoretical analysis has been provided to guarantee the exactness of the relaxation [308]. It is worth pointing out that the result obtained in [308] depended on the popular drag polar, and the approach may need to be re-derived when other drag models are used. The SOCP-based SCP approach has been applied to update the proportional navigation gain for optimal planar engagement with a stationary target by solving a nonconvex OCP online in a receding-horizon fashion with bounded look angle and lateral acceleration as well as impact angle constraint [309].

Later, the midcourse phase of an air-to-ground missile was optimized for maximum impact velocity while locking the target within the missile's field-of-view by an SOCP-based SCP method through combining linearization and convex relaxation with a small-angle assumption [310]. With the aid of multiple techniques such as pseudospectral methods [313, 314, 315], virtual controls [313, 314], and penalty methods [312, 317], the SCP approach has been augmented for solving a wide range of MPG problems including ballistic missile guidance under power system fault and nonconvex thrust magnitude constraint [313], optimal time-varying proportional navigation guidance with impact angle and impact time constraints [311], online midcourse guidance for boost phase interception subject to midcourse-to-terminal handover constraints and heat rate constraint [314], and multi-stage/multi-phase trajectory optimization for dual-pulse missiles with discrete thrust profiles [318, 315]. In particular, to improve the efficiency and robustness of solving optimal MPG problems, convex optimization has been synthesized with MPC to develop a so-called model predictive convex programming (MPCP) method for a class of constrained OCPs in [316] by relating state increments to input corrections and casting the problem as an SOCP problem subject to sensitivity relations. Impact-angle-constrained guidance problems for air-to-ground missiles have been solved as possible applications. In the presence of disturbances and uncertainties, a 3-D interception problem with impact-angle constraints has been addressed by using an intrusive polynomial chaos expansion to transform the stochastic state and constraints into deterministic versions, which has been solved by an *hp*-pseudospectral SCP through combination with a penalty function and backtracking search [317].

### *4.3. Launch/Ascent Vehicles*

The aerospace sector has shown a significant interest in novel launch technologies for safe, efficient, and sustainable access to space. The optimization of the mission trajectories for launch/ascent vehicles (LAVs) is of crucial importance to achieve this goal. However, LAV trajectory optimization is a complex problem due to the highly nonlinear dynamics and stringent mission constraints. Also, the launch trajectory has been usually split into multiple propelled and coasting phases, including vertical ascent, pitchover, gravity turn, fairing jettison, stage separation, coasting, and orbit injection, and proper linkage conditions must be imposed at the transition of each phase (see Figure 11). As such, LAV trajectory design is essentially a multi-phase problem that involves both continuous state variables at the boundaries of the phases and discrete state variables such as the mass of the vehicle due to fairing jettison and stage separations. Further, robust and resilient methods are of paramount importance for LAVs to plan optimal trajectories in real-time onboard to inject the payload into the target orbit with guaranteed accuracy, even under off-nominal conditions due to the presence of model uncertainties and external disturbances such as engine faults or failures. Convex optimization has gained increasing popularity in recent years in addressing these challenges (see Table 10).

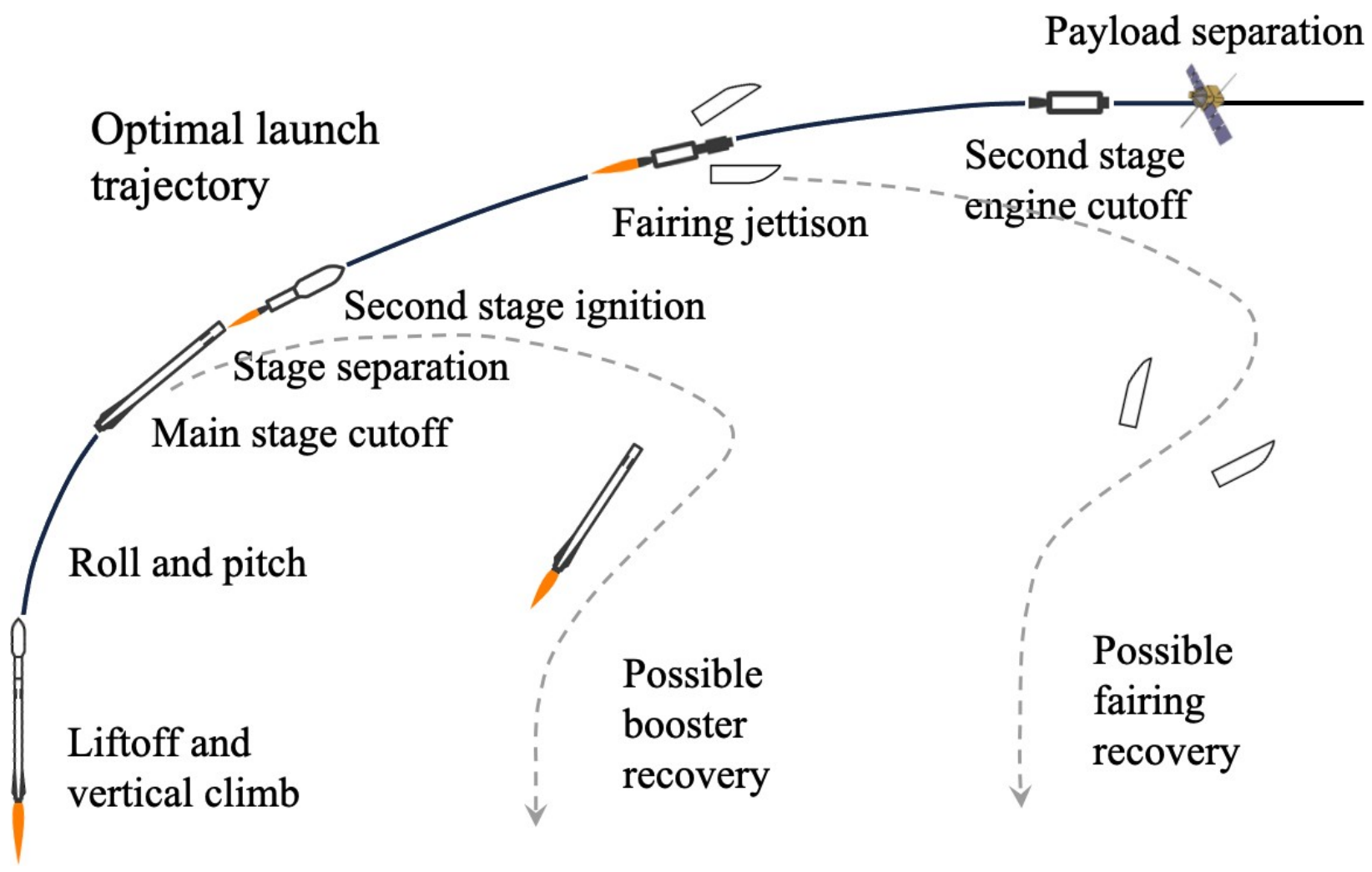

Figure 11: Schematic representation of two-stage launch/ascent vehicles (LAVs).

Table 10: Summary of representative publications on convex optimization for launch/ascent vehicles (LAVs).

| Reference | LAV Problem | Approach | Slack Variables | Change of Variables | Convex Relaxation | Linear Approximation | Solver |
|---|---|---|---|---|---|---|---|
| [319] | 3-DoF maximum terminal velocity and minimum fuel consumption | QP + SCP + pseudospectral | | ✓ | ✓ | ✓ | CVX |
| [320] | 2-D maximum terminal velocity | SOCP + SCP | | ✓ | ✓ | ✓ | SeDuMi |
| [321] | 3-DoF maximum terminal velocity | SOCP + SCP + virtual controls | | ✓ | ✓ | ✓ | MOSEK |
| [322] | 3-DoF maximum terminal velocity | SOCP and LP + SCP | | ✓ | ✓ | ✓ | ECOS |
| [323] | 3-DoF minimum time | SOCP + SCP | ✓ | | ✓ | ✓ | MOSEK |
| [324] | 3-DoF maximum terminal velocity | SOCP + SCP | ✓ | ✓ | ✓ | ✓ | ECOS |
| [29] | 3-DoF multi-phase maximum final mass | SOCP + SCP + virtual controls | | ✓ | ✓ | ✓ | Gurobi |
| [30] | 3-DoF multi-phase maximum final mass | SOCP + SCP + virtual controls + pseudospectral | | ✓ | ✓ | ✓ | Gurobi |
| [325] | 3-DoF multi-phase maximum final mass | SDP + SCP + covariance control | | ✓ | ✓ | ✓ | Gurobi |
| [188, 326] | 3-DoF multi-phase minimum terminal error | SDP + SCP | | ✓ | ✓ | | MOSEK |
| [327, 32] | 3-DoF multi-phase multiple objectives | SOCP + SCP | | ✓ | ✓ | ✓ | MOSEK, ECOS |
| [31] | 3-DoF multi-phase minimum time | SOCP + SCP | ✓ | | ✓ | ✓ | MOSEK |

The first attempt to solve LAV problems using convex optimization appeared in [19], where a successive SOCP approach was proposed with rigorously proved convergence for a type of nonconvex OCPs subject to concave state inequality constraints and nonlinear terminal equality constraints, and an optimal LAV trajectory optimization problem for the upper stage of a medium-lift launch vehicle in a vacuum environment was solved as an application to validate the developed approach. Similar vacuum LAV problems have been solved in [328] as an example to demonstrate the effectiveness of an SDP-based successive convex optimization (i.e., iterative rank minimization) method with guaranteed convergence for a nonconvex QCQP problem transformed from a polynomial OCP formulation of the original problem, in [329] as an application to validate an iterative convex optimization approach with proved convergence based on a Newton-Kantorovich method, in [330] via an iterative SOCP approach by formulating and solving a two-point boundary value problem (TPBVP) through combination of successive linearization and a flipped Radau pseudospectral method, in [331] and [332] using SOCP-based SCP for online optimal LAV trajectory generation in the event of engine faults or failures, and in [333] using an SOCP-based guidance scheme for Mars Ascent Vehicles (MAVs).

In the presence of aerodynamic forces and possible path constraints, more complicated LAV models have been formulated and solved using convex optimization. For example, the iterative convex optimization approach in [329] has been extended to address LAVs under aerodynamic controls subject to constraints on dynamic pressure, axial thrust acceleration, and bending moment through combination of the Newton-Kantorovich method and a Gauss pseudospectral method [319]. Later, an SOCP-based SCP method has been developed to solve a maximum-terminal-velocity LPV problem by approximating the thrust terms as linear functions of the angle of attack and transforming the nonlinear drag coefficient into a linear function of new controls [320]. Moreover, the SCvx algorithm in [69, 70] has been employed to solve LAV problems by modifying the aerodynamic coefficients, introducing new control variables, and relaxing the resulting nonconvex control constraints to facilitate an SOCP formulation of the problem, and virtual controls have also been applied to enhance convergence of the SCP algorithm [321]. Recently, continued efforts have been made to solve LAV problems using convex optimization, including an SCP approach for LAV trajectory optimization by solving a sequence of SOCP and LP problems with as much nonlinearity of the original problem preserved as possible [322], an SCP method enhanced by a modified Chebyshev-Picard iteration discretization technique for minimum-time LAV trajectory optimization [323], and an SCP scheme based on concave-convex decomposition and an augmented Lagrange multiplier method for maximum-terminal-speed LAV trajectory optimization [324].

Due to the multi-stage nature of LAVs, one particular approach to such problems is to formulate and solve multi-phase OCPs, which have also been addressed by convex optimization. For example, by dividing the launch/ascent mission into several arcs and enforcing linkage conditions at the internal boundaries, multi-stage LAV problems can be formulated as multi-phase OCPs, which have been solved by SOCP-based SCP through convenient changes of variables, exact constraint relaxation, and successive linearization with the aid of virtual controls and adaptive trust regions. The approach was initially developed and used to solve a 2-D LAV problem [334] and later extended to address 3-D LAV problems for the SpaceX Falcon 9 launch vehicle [29]. Shortly after, the approach was combined with an *hp* pseudospectral discretization scheme in solving a multi-stage ascent trajectory optimization problem for a VEGA-like launch vehicle subject to nonconvex constraints on the maximum heat flux after fairing jettisoning and the splash-down of the burned-out stages [30]. Also, the approach has been embedded into the MPC framework [335, 336] and the covariance control scheme [325] to gain more robustness to external disturbances and model uncertainties due to engine performance and unpredictable atmospheric conditions. More recently, the full trajectory of a reusable two-stage LAV, including the recovery descent and soft landing of its first stage, has been optimized by this SCP approach [337]. Other recent efforts in using SCP for multi-phase LAVs include QCQP-based SDP relaxation methods for minimum-terminal-error multi-stage launch vehicle trajectory optimization problems [188, 326], SOCP-based SCP methods for trajectory replanning of multi-stage LAVs under dynamic faults such as thrust drop and mass flow loss [327, 32], and an SCP schemed inherited from the Chebyshev-Picard-based SCP [323] for solving a minimum-ascent-time multi-phase

LAV problem [31].

### 4.4. Low-Speed Air Vehicles

Development and deployment of advanced low-speed air vehicles (LSAVs), including UAVs, has gained unprecedented interest in the past two decades for both military and civilian applications. Despite the critical need for mission capabilities such as autonomous operations and online decision-making, most LSAVs are either controlled by onboard/remote pilots or programmed to follow a set of predetermined waypoints. It remains challenging yet highly demanding to develop real-time mission/path planning and trajectory optimization methods as well as resilient G&C strategies to enable optimal, robust LSAV maneuvers and operations for both single-vehicle and multi-vehicle missions, especially in complex, uncertain, and dynamic environments (see Figure 12). In earlier years, MILP and mixed-integer quadratic programming (MIQP) have been used to solve UAV trajectory optimization problems subject to constraints on obstacle/collision avoidance and no-fly zones as well as approximate (usually linear) vehicle dynamics with an aim to enhance its capability of real-time applications [33, 34, 338]. The review in this subsection focuses on the application of convex optimization for LSAV trajectory optimization (not path planning) problems that account for vehicle dynamics and constraints under various representative missions (see Table 11). Particular emphasis is placed on SCP approaches where a series of convex subproblems has to be formulated and solved to find approximate optimal solutions.

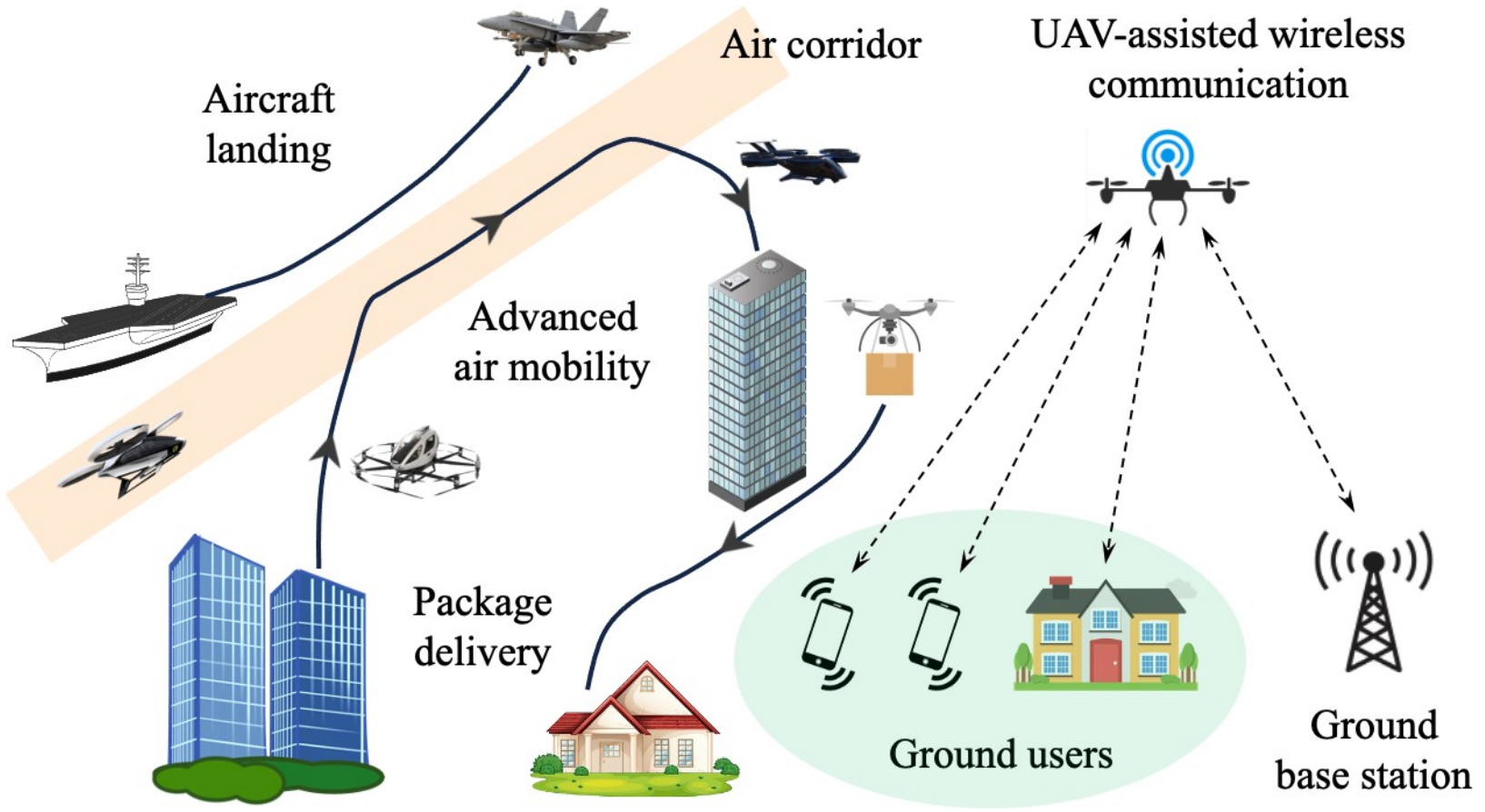

Figure 12: Schematic representation of example low-speed air vehicle (LSAV) missions.

Table 11: Summary of representative publications on convex optimization for low-speed air vehicles (LSAVs).

| Reference | LSAV Problem | Approach | Slack Variables | Change of Variables | Convex Relaxation | Linear Approximation | Solver |
|---|---|---|---|---|---|---|---|
| [339] | 3-DoF minimum total thrust | QP + SCP | | | | ✓ | CPLEX |
| [340] | 3-DoF minimum total thrust | QP + SCP | | | | ✓ | MOSEK |
| [35, 341] | 3-DoF minimum energy | SOCP + SCP | ✓ | | ✓ | ✓ | Bsocp[124] |
| [342] | 3-DoF minimum energy with state-triggered constraints | SOCP + SCP + virtual controls | ✓ | | ✓ | ✓ | ECOS |
| [343] | 6-DoF maximum range and maximum altitude | QCQP + SDP + SCP | ✓ | | ✓ | | SeDuMi |
| [37] | 3-DoF optimal-tracking rendezvous | SOCP + SCP + line search | ✓ | ✓ | ✓ | ✓ | ECOS |
| [344] | 2-D chance-constrained optimal-tracking landing | SOCP + MPC | | | ✓ | ✓ | - |
| [345, 346] | 3-DoF chance-constrained multi-objective | SOCP/MICP + SCP + pseudospectral | ✓ | | ✓ | ✓ | ECOS |
| [36] | 2-D maximum-energy-efficiency communication | SCP | ✓ | ✓ | ✓ | ✓ | - |
| [347] | 2-D maximum-minimum-average-rate communication | SCP | ✓ | | ✓ | ✓ | CVX |
| [348] | 2-D minimum-control-effort AAM | SCP | | ✓ | ✓ | ✓ | MOSEK |
| [349, 38] | 2-D minimum-control-effort multi-phase AAM | SCP | | ✓ | ✓ | ✓ | MOSEK |

The first publication found in using SCP for LSAV trajectory optimization was [339] where 3-DoF collision-free trajectories were generated for a group of LSAVs to transition from a set of initial states to a set of final states satisfying position, velocity, acceleration, and jerk constraints while maintaining a minimum distance between vehicles. The problem was cast as a nonconvex OCP subject to linear dynamics and solved using QP-based SCP by approximating the only nonconvex constraints (i.e., collision avoidance) via successive linearization. However, successive linearization of collision avoidance constraints may lead to infeasible QP subproblems, especially in nonconvex environments. To address this issue, more relaxed, feasible QP subproblems have been formulated, and a decoupled SCP method has been developed for multi-LSAV trajectory optimization by incrementally tightening the collision constraints [340]. Moreover, by introducing a slack variable to relax the nonconvex control constraint, SOCP-based SCP has been developed and demonstrated for 3-DoF quadrotor maneuvering problems [35, 341] and was later extended to address scenarios with compound state-triggered constraints [342]. In addition, an SDP relaxation method similar to that in [244] has been developed for 6-DoF aircraft trajectory optimization [343] and 2-D UAV flight with avoidance zones [350] through a general nonconvex QCQP formulation of the problem. The lower bound of the problem's optimal objective value was sought by solving a transformed SDP problem with a rank-one matrix constraint via an iterative rank minimization approach.

Through common techniques such as change of variables, lossless convexification, convex relaxation, linear approximation, pseudospectral method, and small-angle assumption, convex optimization and SCP have been applied to address various LSAV problems including multi-vehicle formations in centralized [351] and distributed [352] manners with linear dynamics, minimum-time multi-vehicle coordination with nonlinear dynamics and nonconvex obstacle avoidance and inter-vehicle collision avoidance constraints [353], formation rendezvous trajectory optimization of multiple vehicles [354, 355], minimize-terminal-error aircraft landing trajectory optimization [316], stochastic MPC-based aircraft landing under uncertainties and disturbances [344], optimal coordination and rendezvous of unmanned aerial and ground vehicles [356, 37], chance-constrained trajectory optimization for fixed-wing UAVs under probabilistic control and collision avoidance constraints [345, 346], UAV trajectory optimization with avoidance-related constraints via finite-step iteration-free convex reduction techniques [357, 358], and energy management of hybrid aerial vehicles [359].

Among the extensive applications, LSAVs, especially UAVs equipped with efficient convex optimization algorithms, have received much attention in the area of wireless communications by providing cost-effective, flexible, on-demand wireless services such as coverage, relaying, data transmission/collection, wireless sensor network, and internet of things [360, 361, 362, 363]. The key challenge in UAV-assisted communication missions is the optimal balance between maximizing the communication performance (e.g., total information bits transmitted) and minimizing the operational cost (e.g., flight time and energy consumption due to limited battery capacity) while guaranteeing quality-of-service (QoS) and respecting possible constraints on vehicle location, speed, acceleration, and collision avoidance [364, 365, 366]. Building on an energy consumption model, energy-efficient UAV communication problems have been formulated as nonconvex OCPs subject to simple linear dynamics [36, 367]. Through convexification of the objective function and convex relaxation of the nonconvex minimum-speed constraint, an SCP algorithm has been devised to optimize the vehicle's trajectory by jointly considering communication throughput and propulsion energy consumption. Combining with a block coordinate descent method, this SCP approach has been extended from cases with single UAV and single ground user in [36] for multi-UAV wireless systems to serve a group of ground users with maximized minimum throughput [347].

Following similar approaches, the SCP method has been applied for a variety of LSAV-assisted wireless communication problems including UAV-enabled wireless sensor networks by jointly optimizing the vehicle's trajectory and the number of transmission bits [368], minimum-maximum-outage-probability relaying links by jointly optimizing the UAV's altitude, power control, and bandwidth allocation [369], maximum-downlink-sum-rate multi-UAV cellular networks by jointly optimizing resource allocation and base station placement [370], maximum-throughput

UAV-enabled emergency networks by jointly optimizing the UAV's location, power, and bandwidth allocation under statistical QoS constraints [371], UAV-enabled internet of vehicles for intelligent ground transportation by jointly optimizing vehicle communication scheduling, transmit power allocation, and UAV trajectory [372], and secure communication in dual-UAV edge computing systems [373].

In recent years, convex optimization has also found applications for the emerging advanced air mobility (AAM) concept of operations. Enabled by recent advances in battery storage, distributed propulsion, and short/vertical take-off and landing aircraft, AAM aims to explore the third dimension of the space (i.e., airspace) to provide more efficient passenger and cargo air services through urban air mobility (UAM) inside city limits (0-20 miles), sub-urban air mobility (sUAM) connecting a city and its surrounding areas (20-50 miles), and regional air mobility (RAM) for city-to-city transport (50-300 miles). To facilitate fast generation of optimal trajectories for AAM missions, convexification techniques have been introduced to better enable real-time AAM trajectory optimization with initial focuses on single-phase AAM missions [348, 374, 375]. Specifically, through change of variables, convex relaxation, and successive linear approximation, a 2-D minimum-control-effort AAM problem with required time of arrival has been effectively addressed by SCP [348]. Later, the approach has been extended to address multi-phase AAM missions that involve cruise, descent, and landing stages under various operational constraints [349]. Similar to other relevant problems such as LAV missions, the key challenge facing multi-phase AAM lies in the linkage constraints that must be enforced to ensure smooth transitions between phases, which adds another level of complexity to the problem. SCP has shown promising performance in addressing these challenges [349, 38]. In addition, coordinated merging control of multiple AAM vehicles with collision avoidance constraints has also been recently solved by the SCP approach [376]. Studies on convex-optimization-based technique for AAM operations and many other LSAV missions are expected to continue in the coming years.

## 5. Applications to Ground Vehicles

In addition to the aerospace domain, the use of convex optimization for G&C techniques has spread to other vehicular applications such as ground vehicles and intelligent ground transportation for both urban roads and freeways [377, 378, 379, 380]. Many technologies (e.g., sensor, communications, human-machine interface) must work together to enable a safe, efficient, sustainable ground mobility system. In this paper, we focus on motion/speed control and powertrain control, which are critical components of the ground vehicle control (GVC) architecture. The GVC system is expected to take the sensor/navigation data as the input to reapidly generate smooth trajectories, collision-free maneuvers, and optimal control commands for the actuators and powertrains to operate under varying conditions in a dynamic environment such as travelling through an intersection or merging into a main road at minimum energy consumption or control effort with guaranteed safety (see Figure 13). However, these problems are generally nonconvex in their original settings due to the nonlinearity of the vehicle dynamics and the nonconvexity of state and control constraints, which make them difficult to solve in real-time [381, 382]. In this section, we will survey the new schemes enabled by convex optimization that optimize speed trajectories, power split strategies, or both simultaneously for a vehicle in response to the dynamically changing ground traffic environment (see Table 12).

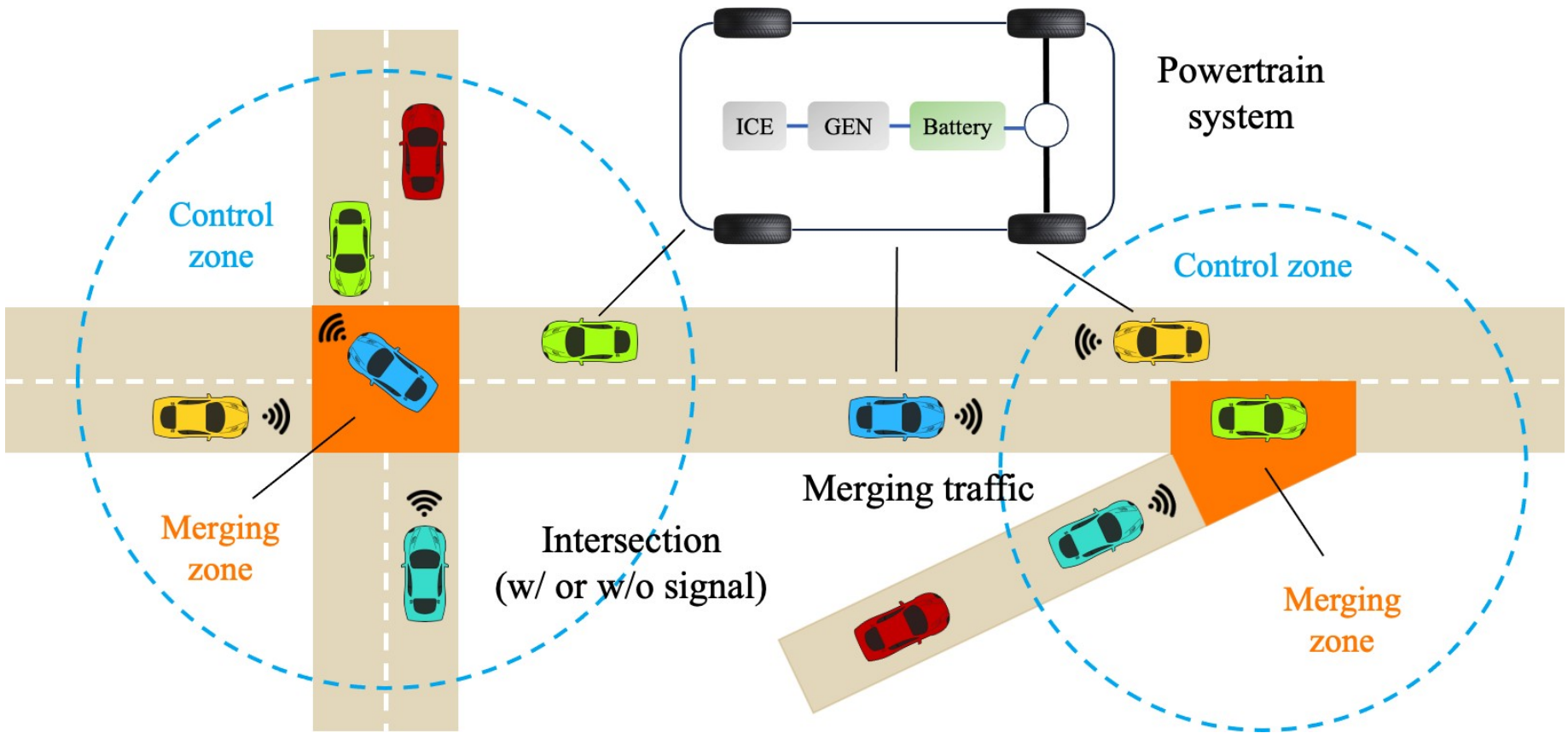

Figure 13: Schematic representation of ground vehicle control (GVC) missions.

Table 12: Summary of representative publications on convex optimization for ground vehicle control (GVC).

| Reference | GVC Problem | Approach | Slack Variables | Change of Variables | Convex Relaxation | Linear/Quadratic Approximation | Solver |
|---|---|---|---|---|---|---|---|
| [383] | Multi-objective speed control | SOCP | ✓ | ✓ | | ✓ | Gurobi |
| [384] | Minimum-time trajectory optimization | QCQP | ✓ | ✓ | | ✓ | FORCESPRO |
| [381, 385] | Optimal trajectory tracking | QP + MPC | | | | ✓ | quadprog |
| [386] | Multi-objective platoon | QP + MPC | ✓ | | ✓ | ✓ | - |
| [387] | Optimal speed tracking platoon | QP + MPC | | ✓ | | ✓ | MOSEK |
| [388] | Optimal trajectory tracking | QCQP + SCP + MPC | | | | ✓ | - |
| [389] | Optimal trajectory tracking | SDP + ADMM | ✓ | | ✓ | | - |
| [40] | Optimal speed control at signalized intersections | QP + SCP + line search + trust region | | | | ✓ | Gurobi |
| [390, 41] | Optimal merging control | QP + SCP + line search + trust region | | | | ✓ | Gurobi |
| [391, 392] | Optimal speed control at unsignalized intersections | QP/SOCP | | ✓ | ✓ | ✓ | GPOPS-II MOSEK |
| [42, 393] | Optimal sizing and control of hybrid electric powertrains | SDP | | ✓ | ✓ | ✓ | SeDuMi |
| [382] | Optimal power-split control | SOCP | ✓ | ✓ | ✓ | ✓ | ECOS |
| [394, 395] | Integrated velocity planning and energy management | QP/SOCP | ✓ | | | ✓ | MOSEK SeDuMi |

When the GVC problem falls into a simple form (e.g., with linear dynamics) that can be easily convexified, a single convex problem or a finite number of convex problems can be solved to get the solution. For example, a collision avoidance problem has been solved to generate the optimal trajectory and the corresponding force and moment to be distributed to each tire via an SOCP approach based on an assumption of non-rotating point-mass vehicle [396]. This assumption has been relaxed by considering the yaw of the vehicle, which adds significant nonconvexity and complexity to the collision avoidance problem. Instead of solving a single convex optimization problem to get the solution, a three-stage approach has been proposed by solving a convex optimization problem in each stage through a series of convex approximations [397]. In addition, SOCP has also been used to generate optimal speed profiles for vehicles moving along a fixed trajectory subject to affine dynamics, hard constraints on friction circle, speed limit, and time window, and semi-hard constraints on comfort acceleration and deceleration in both static and dynamic environments [383]. Through slack variables and penalty functions, single SOCP problems have been solved with preserved convexity and global optimality, balancing multiple performance metrics including smoothness, time efficiency, and speed deviation. In emergency situations, vehicle trajectory and the associated control inputs need to be replanned in real-time at possible friction limits to achieve the minimum response time. Assuming the existence of a nominal trajectory, a minimum-time OCP has been solved considering road topography as well as engine power and tire friction limits to replan the vehicle's speed and trajectory in emergency obstacle avoidance scenarios by approximating the problem as a convex QCQP based on a simplified point-mass vehicle model [384].

When it comes to more complicated GVC problems subject to nonlinear vehicle dynamics and nonconvex constraints such collision avoidance, successive linearization has been frequently used to convexify the problems into more favorable forms that can be potentially implemented in real-time online under MPC frameworks by solving a convex problem in each circle. For example, tailored MPC algorithms have been developed through convex QP approximations to the original nonconvex problems for optimal trajectory tracking by controlling the front steering angle [385] and the braking torques at the four wheels subject to collision-avoidance constraints and physical limitations on the actuators [381]. These techniques have been incorporated into shared control schemes where the vehicle is commanded by the controller and a human driver in a safe manner [398]. In addition to environmental obstacles and the handling limits of the vehicle, the driver's intent can also be integrated into the MPC framework where a set of convex problems is iteratively solved [399]. Due to the quadratic nature of the collision avoidance constraints and quadratic objective functions frequently used in MPC, dynamic collision avoidance problems for a group of vehicles have been formulated as a nonconvex QCQP, which has been solved by SCP through affine approximation of the nonconvex collision avoidance constraints to produce optimal commands for the vehicles to follow a reference trajectory while avoiding collisions between the vehicles [388].

In the context of connected and automated vehicles (CAVs) enabled by vehicle-to-vehicle (V2V) and vehicle-to-infrastructure (V2I) communication technologies, convex optimization has been used for vehicle motion control under various scenarios such as platoon or car following. For example, by modeling the platoon as a multi-agent interconnected dynamic system with shared traffic information, [39] solved a convex QCQP problem with double-integrator longitudinal dynamics to generate the optimal maneuver of each vehicle in a decentralized manner subject to speed, control, and safe-distance constraints. Using an optimal velocity with relative velocity (OVRV) car-following model, maneuvers of CAVs have been optimized in mixed-autonomy multi-lane traffic scenarios considering traffic efficiency and driving comfort of both CAVs and human-driven vehicles (HDVs) [386]. The problem was initially formulated as a nonconvex mixed-integer programming problem, where the nonlinearity and nonconvexity come from the car-following dynamics of the HDVs, and the integer variables stem from the lane-change decision variables. The problem was then transformed into a convex QP problem through a linear approximation of the piecewise linear car-following dynamics and relaxation of the integer lane change variables into three separate convex subproblems [386]. Recently, a distributed tube-based MPC method has been developed for platoon control of heterogeneous CAVs in the presence of modeling uncertainties and measurement disturbances [387]. Through proper relaxation of the nonlinear vehicle dynamics into linear forms via a change of independent variable and state transformations, convex QP problems have been established and solved locally by each CAV for optimal speed tracking maneuvers.

Bottleneck GVC challenges at signalized intersections and merging roadways in the context of CAVs have also been addressed by convex optimization to improve traffic efficiency while ensuring safety. Signal phase and timing (SPaT) information from the upcoming traffic lights enables predictive planning and control of CAVs to pass the signalized intersections safely and efficiently with reduced fuel consumption and travel time. However, it still remains a solid challenge to generate optimal speed and control profiles for the vehicles moving along signalized corridors obeying dynamically changing SPaT and satisfying other constraints such as speed limit, control limit, and collision avoidance. A mixed-integer SCP method has been developed for optimal speed control of CAVs over multiple signalized intersections, and the integer variables stem from the selection of green phase window to cross the intersection [400]. However, mixed-integer programming problems are NP hard, and the computational speed significantly drops when the number of integer variables and the problem size grow. To remove the integer variables, the green window selection problem has been reduced to the determination of the reference velocity from the upcoming SPaT and the distance between the vehicle and the intersection [401]. Based on this strategy, SCP algorithms enhanced by line search and trust region techniques have been developed to minimize the fuel consumption while avoiding idling and frequent stop-and-go patterns of CAVs at signalized intersections considering complex inter-vehicle interactions and nonlinear vehicle dynamics [40]. With the aid of pseudospectral discretization, the computational efficiency of the designed SCP algorithms has been greatly improved, enabling real-time on-vehicle applications under MPC frameworks for resilient response to emergency situations. In combination with a rule-based merging sequencing strategy, these improved SCP algorithms have recently been extended to determine the optimal speed profiles of CAVs on merging roadways [390, 41].

Movements of CAVs at unsignalized intersections can also be coordinated for optimal traffic mobility. For example, a nonconvex OCP problem has been formulated and solved in [391] by convex optimization for multiple CAVs optimally and cooperatively crossing a signal-free intersection. Convexification techniques have been used to transform the problem into a convex QP problem while remaining as much nonlinearity as possible. Specifically, the nonlinear dynamics are linearized through a change of independent variable (from time to distance travelled), domain transformation (from time to space), and a change of variable (from velocity to kinetic energy). The nonlinear relationship between kinetic energy and velocity is approximated by a linear function. The conservativeness of the convexified QP formulation has been verified through comparison with the original nonconvex problem [391]. More recently, the CAV coordination problem at unsignalized intersections has been addressed through a two-level hierarchical approach where the upper level determines an optimal crossing order while the lower level optimizes the speed trajectories of all CAVs with guaranteed collision avoidance following the crossing order from the upper level [392]. Based on the domain transformation approach and the change of variable in [391], both the upper- and lower-level nonconvex OCPs have been relaxed into SOCP problems that seek trade-offs between energy consumption and travel time. In addition to linear approximation and domain transformation, semi-definite relaxation has also been applied for cooperative planning and control of multiple CAVs in unsignalized multi-way junction and intersection scenarios. For example, a nonconvex GVC problem has been formulated and divided into two small subproblems subject to nonlinear dynamics and nonconvex coupled collision-avoidance constraints, respectively, using the ADMM method [389]. The subproblem with inequality nonconvex collision-avoidance constraints was relaxed into an SDP problem. The optimality and feasibility of the solution to the original nonconvex problem may not be guaranteed by the solution of the relaxed SDP problem, although the SDP method can usually provide accurate or near-optimal approximations with higher computational efficiency than other nonconvex approaches [389].

The way how the vehicle operates has obvious effects on the fuel and/or electric consumption that can also be optimized through convexification approaches by re-placing the nonconvex feasible set with a convex superset, enabling the solution to the original problem by efficiently solving a relaxed convex problem [43, 382]. A typical example of such relaxations can be found in [402, 42], where a systematic convexification approach has shown to be promising to efficiently solve a highly coupled, nonconvex, mixed-integer problem of simultaneously optimizing battery size and energy management for a plug-in hybrid electric vehicle (HEV). Through a change of variables for equivalent convex relaxation of the battery model and a convex second-order approximation of the engine-generator unit, the original nonconvex problem was transformed into a convex problem, which was solved in each iteration of two nested loops (one through all given sizes of engine-generator unit and electric machine and the other through all possible distributions of charging stations along a known bus line) to obtain a solution near the global optimum [402, 42]. The approach has been extended to devise a heuristic method for simultaneous optimization of battery dimensioning and power split of a plug-in HEV by first deciding the feasible values of the integer variables (engine on/off control) and then solving a convex subproblem to obtain the optimal values of the remaining design variables [393]. The strategy has been shown to be able to converge toward a solution of guaranteed optimality using Pontryagin's minimum principle [403] with much higher computational efficiency than dynamic programming [42, 393]. Further, the approach has been employed to optimize both the powertrain size and power management of fuel cell HEVs [404] with models of different levels of details [405] by formulating and solving SOCP problems. More recently, with the aim to significantly reduce the computational burden for real-time applications, convexification methods have been used for optimal energy management of power-split HEVs [406], optimal power allocation of HEVs in combination with ADMM [407, 408], energy management of HEVs with battery degradation through SOCP and MPC [409], and integrated speed planning and energy management of autonomous fuel cell HEVs [395] and connected fuel cell HEVs passing through signalized intersections [394]. The role and application of convex optimization for component sizing and energy management of HEVs have been discussed in a recent review paper, and interested readers are referred to [410] for details.

## 6. Future Research Directions

In this section, we will discuss some issues, challenges, and future research directions related to the application of convex optimization for vehicular G&C problems.

### *6.1. Theoretical Advancement of Convexification Techniques*

The existing literature on SCP-type methods mainly focuses on using numerical simulations to validate their real-time performance, optimality, and convergence in solving various vehicular G&C problems. However, theoretical development and convergence proof of the SCP method receives very limited attention. While a few attempts have been made to solve nonconvex OCPs through successive convexification with guaranteed convergence properties, the theoretical guarantees are usually based on the assumptions of special dynamics (e.g., linear or control-affine systems) and constraints (e.g., concave state inequality constraints) [19, 71, 74]. It is challenging yet valuable to explore more advanced SCP algorithms with theoretically proved convergence by relaxing assumptions on the problem settings while maintaining lossless convexification to expand the class of problems that can be handled [65]. A key next step is to develop enhanced versions of the SCP algorithm by leveraging other techniques and safe-guarding mechanisms such as virtual controls, virtual buffer zones, line search, and trust region to construct more comprehensive SCP frameworks for more general OCPs with nonlinear dynamics and nonconvex state and control constraints to be relaxed and solved by polynomial-time convex optimization methods [70]. Thorough analysis of the convergence properties (e.g., weak or strong convergence and convergence rate) of these SCP algorithms are expected, and numerical simulations of nonconvex example problems are then needed to validate these theoretical results. In addition to guaranteed convergence, future work will also need to focus on the proof of the exactness of the utilized lossless convexification or convex relaxation techniques from theoretical perspectives by showing that the relaxed problem is equivalent to and share the same solution with the original problem. Moreover, other fundamental issues, including the feasibility of each subproblem parameterized and solved within SCP, the effects of the feasibility of subproblems on the convergence of SCP, the existence of optimal solutions to subproblems, and the quantification of time and space complexity of the problem, are also valuable to be explored to gain more certainty, transparency, and confidence in the performance of the algorithm.

### *6.2. Fundamental Improvement of Convexification Techniques*

While convex optimization and SCP algorithms have gained significant popularity as effective methods for solving a number of vehicular problems from different domains, fundamental issues exist and need to be addressed to further improve their performance such as convergence, robustness, and accuracy to enable more reliable and efficient vehicular operations in uncertain, dynamic mission environments. One of the biggest challenges is to provide a good initial guess for the SCP method. No user-provided guesses are needed for convex optimization algorithms such as IPMs to solve a single convex problem; for SCP approaches, however, good initial guesses are required, and the convergence and results greatly depend on the initial guess. Perhaps the best initial guess is the actual optimal solution, however, the solution of the problem is not known *a priori*. Therefore, techniques and strategies are desired to design suitable initial guesses for the SCP process with better convergence. For example, the convergence of SCP may be accelerated by infusing with the indirect methods [206]. The indirect methods convert the original OCP into a two-point boundary value problem by formulating the necessary conditions for optimality based on Pontryagin's minimum principle and may converge very quickly to the optimal solution. However, the indirect methods are very sensitive to the initial guess. The initialization of the indirect methods may be mitigated by extracting information from the multipliers at each SCP iteration, and the resulting indirect methods may then be combined with SCP to decrease the total number of iterations required for convergence [74]. Other strategies such as continuation or homotopy may also be effective in bypassing the need for a good initial guess for SCP [193, 210]. In addition, the existing convex optimization approaches greatly rely on linear approximations of nonlinear dynamics and nonconvex constraints for convexification purposes and then apply SCP to solve the problem. However, simple linearization may result in poor approximations of the original nonlinear formulations, making the SCP approach more sensitive to the initial profiles and more difficult to converge [322]. Trust-region constraints have been imposed as routine techniques to improve convergence; however, the trust-region radius is a key parameter that needs to be carefully adjusted. As such, it is rewarding to further explore the structure of the problem and develop new convexification methods that can retain as much nonlinearity of the original problem as possible such that the SCP method can quickly converge in fewer iterations even without good initial guesses or trust region strategies. Furthermore, future work should also focus on combination and implementation of convex optimization methods with MPC or covariance control frameworks to explicitly incorporate uncertainties and disturbances in a closed-loop manner in the design of G&C systems for more robust and resilient vehicular operations with real-time performance.

### *6.3. Customization of Convex Optimization Algorithms*

In this literature, generic solvers, such as SeDuMi, ECOS, MOSEK, Gurobi, and SDPT3, have been used to solve the resulting convex problems. It is fine focusing the effort on the convexification process and using the off-the-shelf solvers to demonstrate the performance of the approaches in the initial stages of the development. When the stability and efficiency of the methods have been validated, however, specific effort need to be made to tailor these algorithms to solve specific application problems for verification and real-world implementation purposes. Autonomous vehicular systems rely significantly on onboard computation and optimization to operate in a range of scenarios and environments safely and efficiently with stringent real-time requirements and limited memory. Therefore, taking full advantage of specific problem structures to develop customized convex optimization algorithms with significantly reduced number of mathematical operations and computational branches would be critical and of great interest in the future for software verification and embedded system applications [124]. For example, sparsity can be explored and leveraged to translate the problem into sparse optimization formulations to minimize memory usage and number of arithmetic operations for increased computational speed [219, 218, 282]. Other techniques such as approximate minimum degree ordering [411] and explicit coding [123] may also be used to increase the computational speed of IPMs in solving the relaxed convex optimization problems.

### *6.4. Extension to Formulations of Higher Fidelity*

The aerospace and automobile industries have an explicit interest in optimizing vehicle control and operation involving multiple disciplinary models such as aerodynamics, vehicle dynamics, structural dynamics, and propulsion. However, simultaneously integrating multi-physical, coupled disciplines in G&C of vehicles is challenging due to their complex dynamical interactions and the computational requirements of high-fidelity models [412, 413]. Due to the lack of efficient means to directly integrate computationally expensive, high-fidelity models in G&C designs, the existing G&C methods for practical real-time applications, specifically for air, surface, and underwater vehicles, have been limited to extensive use of low-fidelity, low-dimensional dynamical models that consist of rigid-body vehicle dynamics with analytical disciplinary models (e.g., aerodynamics), ignoring many important physical processes (e.g., rotor-rotor interaction, complex gust profiles) [414]. The model of reduced-fidelity may result in poor closed-loop performance in control, e.g., suboptimality, instability, constraint violations, and lack of robustness, and such suboptimal performance is unacceptable when it comes to critical applications such as aircraft landing in highly disturbing environments [415]. Therefore, a technical gap remains that pertains to the resolution of the dilemma between model fidelity and computational efficiency. Recent advancement of IMPs and convexification techniques provides significant opportunities to develop innovative, rigorous G&C algorithms that seamlessly integrates accurate yet fast dynamical models to facilitate more reliable and efficient model-based vehicle G&C and decision-makings with guarantees on stability, optimality, computational efficiency, and robust constraint satisfaction in the context of vehicular missions with high-dimensional models. Future research may be directed toward extending the convex optimization and SCP approaches to more complicated problem formulations in more realistic scenarios with higher-fidelity vehicle models. Example applications include hypersonic trajectory optimization with high-fidelity aerothermodynamic models, aircraft landing with strong aerodynamic perturbations and complex dynamical responses, and speed/motion control of ground vehicles considering high-fidelity energy consumption models. All these situations would result in problems of increased complexity, and mindful convexification approaches should be developed to retain the overall computational efficiency of the solution process.

### *6.5. Integration with Other Systems and Techniques*

While convex optimization has been successfully applied to address a number of vehicular G&C problems, it remains worthwhile to integrate these algorithms with other relevant systems (e.g., navigation, sensing, communication) and techniques such as data-driven modeling, machine learning approaches, and even indirect optimal control to further enhance the system performance for fully autonomous vehicular operations in complex mission environments. For example, the SCP method can bypass the initialization issues of the indirect methods, which in turn would promote the convergence of SCP, as discussed in subsection 6.2. In addition, unified, end-to-end guidance, navigation, and control (GNC) systems may be designed by integrating these traditionally isolated disciplines under efficient convex optimization architectures. However, a complex and open question concerns the extent to which these disciplines should be integrated. Recently, there seems to be a growing interest in the combination of convex optimization and machine learning as integral components of G&C loops to enable rich use of online computational methods for improved mission performance and robustness. For instance, SCP has been used to generate optimal trajectories for deep neural networks (DNNs) to learn in multi-agent space motion planning missions [416]; DNNs have been used to predict the optimal flight time [417] or generate the initial guess [418] for SCP algorithms to optimal PDG problems with improved efficiency; DNNs have been designed to map any observed actual flight state of the spacecraft to optimal RPO actions through training on a set of optimal trajectories generated by convex optimization [419];

Convex optimization has been combined with a reinforcement learning framework to design optimal low-thrust lunar transfers [420]; A neural network has been designed to approximate solutions of a centralized method for multi-spacecraft coordination using training data generated by a centralized convex optimization framework [421]; Also using datasets generated convex optimization, DNNs have been developed for mission reconstruction of launch vehicles under thrust drop by mapping from the fault state to the optimal rescue orbit [422]; A hybrid framework has been designed by incorporating reinforcement learning and convex optimization to cooperatively solve a UAV-based data collection problem [423]. In addition, data-driven approaches (e.g., reduced-order modeling) have been investigated in recent years to obtain accurate dynamical models by capturing the previously uncaptured complex factors. It would be of interest to validate the convex optimization approaches on purely data-driven models or mixed systems with both model-based and data-driven terms.

### *6.6. Application to Multi-phase Vehicular Missions*

An additional problem that is under-explored concerns vehicular operations involving multiple mission phases/stages. A typical example is the multistage launch vehicle ascent problem. The ascent trajectory has been divided into multiple flight phases from liftoff to payload release, consisting of a sequence of propelled and coasting arcs and featuring mass discontinuities due to the separation of inert masses at stage burnout [29, 30]. Effective G&C schemes accounting for variable conditions need to be developed to meet all mission requirements. Also, we have discussed planetary entry and powered descent/landing separately in the previous sections. These mission phases are essentially connected in the entry, descent, and landing (EDL) architecture primarily for Mars exploration. It would be of interest to piece these phases together and integrate the algorithm of each portion for more optimal EDL mission designs. Specifically, some Mars missions may have a parachute descent phase before the powered descent phase initiated by parachute cutoff. The parachute cutoff time is a critical parameter that may be optimized in the multi-phase EDL mission framework [77]. Another example is the recent AAM concept of operations. A complete AAM mission profile may involve multiple flight phases including takeoff, ascent, merge, cruise, descent, and landing. When the flight schedule is determined and passed to the vehicle, it is critical for the vehicle to generate and follow an accurate, smooth trajectory across all the possible phases to safely arrive at the destination vertiport or airport at minimum energy consumption while meeting all the stringent time and regulation constraints [349, 38]. For ground vehicles, if the distance to the next intersection is obtained and the SPaT information of the upcoming traffic lights is known, the speed profile of the vehicle can be optimized to best pass the intersection. Comparing to isolated intersections, multi-intersection traffic control at signalized corridors seems more attractive. It is challenging yet beneficial to optimize the maneuvers of the vehicle to pass a series of green traffic lights without having to stop at red traffic signals (also known as green wave) for maximum energy savings and minimum greenhouse gas emissions [40]. All these cases can be potentially handled as multi-phase OCPs, which may be addressed by convex optimization methods. Because the problem consists of multiple phases, one major concern would be to enforce proper linkage conditions at the boundary between adjacent phases.

### *6.7. Application to Wider Vehicular Missions*

In addition to space, air, and ground vehicular missions surveyed in this paper, future effort is expected to push the boundaries of convexification techniques and extend the list of problems that can be addressed by convex optimization methods. In fact, several publications have been found in the maritime domain using convex optimization to address G&C problems for surface and underwater vehicles. For example, a collision-avoidance problem for multiple unmanned surface vehicles (USVs) has been formulated and converted to a convex QP problem that can be solved in real-time for safe and optimal direction and path of each vehicle in the presence of static or moving obstacles [424]. In addition, a problem of jointly optimizing trajectory and communication resource allocation has been addressed for a USV-enabled maritime wireless network where the USV is used to assist the communication between the terrestrial base station and ships [425]. Considering the USV kinetics and multiple constraints such as safe sailing, breakpoint distance, line-of-sight links, and resource allocation, a joint optimization problem that maximizes the minimum throughput among all the ships has been established and decomposed into two subproblems that are solved iteratively using successive convexification and IPMs. For underwater vehicles, a 3-D trajectory tracking problem has been considered based on a 6-DoF dynamical model and transformed into a convex QP problem [426]. To improve the robustness of the tracking control method under model uncertainties and disturbances, the problem was solved within the MPC framework. This QP-based approach has later been implemented in a so-called double closed-loop MPC scheme for underwater vehicle trajectory tracking [427]. The outer-loop position controller generates the desired speed command that is tracked by the control forces and moments produced by the inner-loop speed controller. More recently, convex optimization has been applied to solve a position tracking problem for an underwater vehicle to track the reference trajectory under attitude and velocity constraints in the presence of saturated thrusts and time-varying disturbances [428]. G&C of surface and underwater vehicles considering high-fidelity hydrodynamic models are also worth investigation as discussed in subsection 6.4. More vehicular mission scenarios, such as aircraft taxi [429], G&C of robotic rovers and helicopters on Mars [430], shipboard landing of naval aircraft considering ship motion and coupled ship-aircraft airwake [431], and vehicular operations in the presence of severe disturbances and faults/failures of the actuators or other components [432], may also be studied and addressed by convex optimization techniques in the future.

### *6.8. Application to Cross-domain Vehicular Missions*

The existing research on convex-optimization-based G&C problems has been almost invariably limited to single domains by focusing on space, air, or ground missions separately with very few exceptions. Joint missions across multiple domains using heterogeneous vehicular platforms may bring unprecedented mission performance that may not be achieved by a single type of vehicles. In recent years, there has been a growing trend towards utilization of heterogeneous multi-agent systems for a variety of military and civilian applications such as emergency response, search and rescue, and cooperative communication and sensing. Future research directions may involve studies on application of convex optimization methods for cross-domain vehicular missions. For example, considering the agile maneuverability, broad field of view, and rapid coverage of large areas of UAVs and the accurate location control of unmanned ground vehicles (UGVs), there has been a significant interest in exploiting the complementary capabilities of aerial and ground vehicles to develop cooperative UAV-UGV systems have can be used for different applications such as surveillance, communication, UAV refueling/charging, sidekick package delivery, and target detection and tracking. As a specific application scenario, a UGV can be deployed as a mobile base station to refuel/charge the UAVs at different locations to reduce idle time and increase the utility of fuel/energy. A similar situation is the sidekick UAV-UGV delivery system where the UGV serves as a mobile hub that travels on the main roads, stops at specific locations, and dispatch UAVs to deliver packages to the regions with natural disasters or limited road networks. Both cases may involve a coordinated motion problem where the UAV may rendezvous with and land onto a moving UGV, which has been solved by SCP algorithms based on the error dynamics [356, 37]. In addition, both satellites and UAVs have been employed to assist with ground communication by constructing integrated space-air-ground networks through coordinated transmissions and using UAVs as relays for data transmission between satellites and the ground facility [433, 434, 435]. The system capacity has been maximized by jointly optimizing transmit power allocation, device connection, and UAV trajectory via SCP, and satellite maneuvers may be considered in the future.

### *6.9. Physical Demonstration and Experimentation*

Finally, the existing work on convex-optimization-based G&C has been almost focusing on algorithmic development and validation in simulation environments. Even though the preliminary results on the effectiveness of the developed G&C methods are encouraging, a gap between theoretical development and convincing experimental results has been observed, which drives the necessity of performing exten-sive validation and verification campaigns and real-world experimental tests over a wide range of mission conditions and scenarios to comprehensively demonstrate the capability of the algorithms. To this end, testbeds should be designed to emulate vehicular motion systems while reflecting realistic mission scenarios subject to hardware limitations such as limited on-board memory and computing power, which may put extra constraints on the development of implementation of the algorithms [436, 437]. In addition, custom peripheral devices may be preferred over high-performance general-purpose microprocessors to carry out all or part of the G&C function for experimentation, due to the fact that hardware applied for vehicular applications (especially for space missions) usually lags that for general-purpose usage [438]. Some programmable hardware and devices, such as field programmable gate arrays, may be suitable for customization and prototyping purposes for specialized vehicular applications. In the future, methodologies that are comprehensively validated in simulated environments are expected to be further tested and evaluated through physical experiments on testbeds either in controlled lab environments or on real-world platforms. Key aspects to be assessed include the computational cost of G&C command generation, optimality and accuracy of the produced solutions, and the ability to consistently deal with modeling errors, uncertainties, and external disturbances. The continuous effort in developing more powerful computing hardware with higher update frequencies, along with the ongoing research in convex optimization, paves the way for the use of computational convex-optimization-based G&C algorithms for future vehicular operations.

## 7. Conclusions

This paper provides an overview of convex optimization approaches and surveys their applications to the design of G&C algorithms for space, air, and ground vehicle systems. Convex optimization enables real-time computation of optimal or suboptimal solutions, provides better G&C capabilities, and enhances opportunities for more efficient vehicle operations with improved overall performance. The motivating factors that drive the development of convex optimization techniques for modern G&C systems have been summarized, and the existing challenges and issues that play a central role in the evolution of these approaches have been identified and discussed in this paper. In each vehicular domain surveyed, a wide range of nonconvex G&C problems has been systematically transformed into and solved as convex problems through a series of convexification techniques such as change of variables, convex relaxation, and successive linearization. The main purpose of this paper is to stimulate and promote the interest of G&C researchers to apply their expertise to advance the next-generation G&C technologies using convex optimization. The convex-optimization-based G&C field is still rapidly evolving and may see deeper theoretical advancement, wider applications and physical implementation, and more real-world deployments over the next decade. We expect this paper to encourage discussion regarding the future direction of this area.

## 8. Acknowledgements

The author gratefully acknowledges the support to this work by the National Science Foundation (NSF) grants CNS-2231710, CMMI-CAREER-2237215, and EEC-2244052; the Office of Naval Research (ONR) grant N00014-23-1-2607; the National Aeronautics and Space Administration (NASA) University Leadership Initiative (ULI) grant 80NSSC23M0059; and the Oak Ridge National Laboratory (ORNL) subcontract CW41203.